\documentclass[a4paper,11pt,times]{article}
\usepackage[dvips]{graphics,graphicx,epsfig}
\usepackage{amsfonts}
\usepackage{amsmath}
\usepackage{amssymb}
\usepackage[english]{babel}
\usepackage{graphicx}

\usepackage{amsthm}
\usepackage{times}
\usepackage{boxedminipage}
\usepackage{graphics}

\usepackage{epsfig}
\date{}    %%%% so that date does not print.
\newtheorem{theorem}{Theorem}

\newtheorem{definition}[theorem]{Definition}

\newtheorem{defn}[theorem]{Definition}
\newtheorem{algo}[theorem]{Algorithm}

\begin{document}

\begin{center}
\huge{\textbf{Dynamic bandwidth Reservation in Virtual Private Network under Uncertain Traffic}}\\

\bigskip

\Large{H\'el\`ene Le Cadre}

%\bigskip

%\large{helene.lecadre@enst-bretagne.fr}

\bigskip

ENST Bretagne\\
Technop\^ole de Brest Iroise\\
FRANCE\\
\small{helene.lecadre@enst-bretagne.fr}\\

\end{center}
%%%%%%%%%%%%%%%%%%%%%%%%%%%%%%%%%%%%%%%%%%%%%%%%%%%%%%%%%%%%%%%%%%%%%%

%\maketitle                        %%%% To set Title and Author names.
\thispagestyle{empty}

%%%% Replace with your Abstract.

\noindent
{\bf Abstract -
The aim of this paper is to analyze the dynamic evolution of a Virtual Private Network. The network is modeled as a system, controled by a manager who should take appropriate decisions. However, to be able to take the best possible decisions, the manager should also be able to forecast the worst behavior, in the sense of a quality of service criterion, of the system, he wants to control. We have chosen to model this problem, as an iterative two side game. On the one side, the operator tries to reserve the minimal amounts of bandwidth to guarantee the best possible quality of communication for its various clients. On the other side, the traffic of the clients follows the worst behavior, in the face of the reserved bandwidths.\\
The theory of Markov decision processes (MDP) enables us to model the uncertainty associated to the knowledge of the traffic. Besides, two levels should be differentiated in our system. The local level of the clients, who evolve independently of one another and selfishly, choosing the worst possible traffic evolution. At this level, the manager could reserve bandwidth locally, on each link for every Virtual Private Network. Whereas, at the global level of the links, decisions should be taken by the manager to centrally control the network.\\
A hierarchical MDP approach and the stochastic game framework are introduced to propose solutions to this difficult problem. Furthermore, we study the asymptotic behavior of the system, and prove the convergence towards stationary strategies. In the final section, we introduce parametrized strategies, whose parameters should be estimated with the help of simulation. Indeed, simulation based optimization, over the policy space, provides us an alternative to Bellman's principle, all the more interesting as this principle might become hard to apply, when the cardinality of the state space increases.
  
{\small\em }
}

\vspace{0.5cm}

{\small Hose model, Markov Decision Process, Bellman's optimality principle, stochastic Games, Cross-Entropy method}

%%%%%%%%%%%%%%%%%%%%%%%%%%%%%%%%%%%%%%%%%%%%%%%%%%%%%%%%%%%%%%%%%%%%%%%%%%%%%%%%%%%%%%%%%%%%%%%%%%%%%%%%%%%%%%%%%%%%%%%%%%%%%%%%%%%%%%%%%%%%%%%%%%%%%%%%%%%%%%%%%%%%%%%%%%%%
\section{Introduction}
During the last decades, many methods have been developed to tackle the rather hard problem of traffic matrix estimation. Our purpose in this article is not to develop a new method for traffic matrix estimation, but rather to consider the problem under a system oriented point of view. Indeed, our system is made of a telecommunication network of nodes and directed links. The operator, or the network manager has the possiblity to act on the bandwidth reservation, in view of the evolution of the traffic going through the whole network. We assume that, at each global bandwidth allocation, the traffic evolves, following \textit{the worst configuration} in the sense of a Quality of Service (QoS) criterion. The network operator should be able to forecast the worst possible evolution of the traffic, and to propose solutions so as to drive the network in an optimal way. In the context of Virtual Private Networks, guaranteeing an admissible QoS, via reserved bandwidths, loss, and delay characteristics, is a crucial task for the network manager.\\   
Virtual Private Networks (VPNs) are networks built between geographically distant IP-sites of a firm. With the help of this technology, distant sites of the same firm are able to communicate via secured tunnels. Indeed, the data should be transmitted via Internet, which is a public infrastructure shared by many operators. In order to guarantee the security of its client, the data will be encrypted and sent along virtual tunnels using MPLS technology. Besides, a Service Level Agreement (SLA) contract should be passed between the network provider and its client. The aim of this treaty is to specify bounds on admissible levels of QoS. As a result, the manager should be able to forecast both the spatial and the temporal evolution of its traffic.\\
Traditionaly traffic matrices are used to solve such problems. Nevertheless, their accuracy rely mainly on the quality of the estimator itself and of the data, which can be quite hazardous. The solution we have chosen to get a rough characterization of the traffic, is to use the hose model, introduced for the first time in $[1]$.\\ 
The client is asked to merely specify:\\
-the amount of traffic going in/out each of its web sites,\\
-the relationships between all its web points (source $\rightarrow$ destination).\\

\begin{center}
\begin{figure}[h]
\includegraphics[scale=0.5]{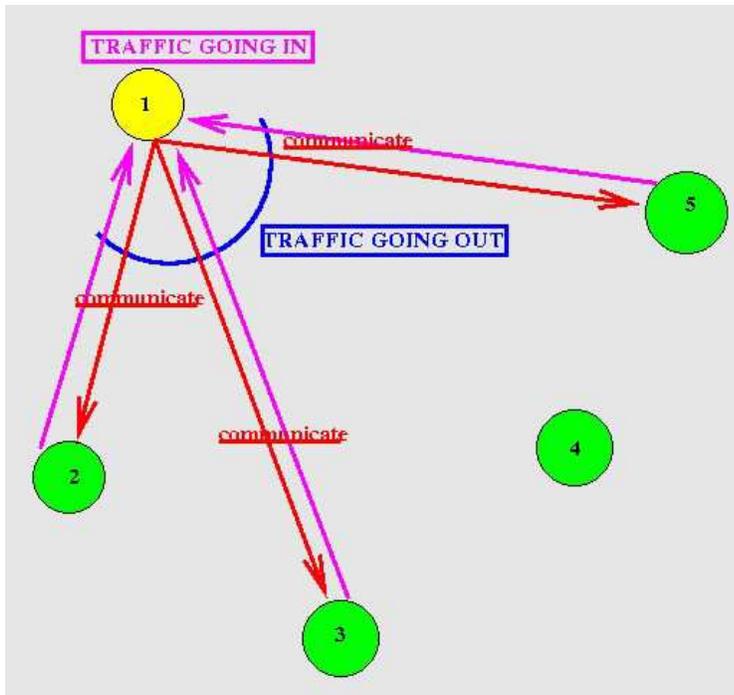}
\caption{The hose model.}
\small{As an example, we consider the hose model applied to a small network. The firm is composed of $5$ different IP sites, which are supposed to be geographically distant.
For the site $1$, which is supposed to be the head of the firm, the client gives the operator the connections to the other areas: $1\leftrightarrow 2,\;1\leftrightarrow 3,\;1\leftrightarrow 5\;,$ where the symbol $\leftrightarrow$, means that there is a potential bidirectionnal connection between the two sites. Besides, the client gives the volume of traffic going out of the site $1$, and, possibly, the amount of traffic going in the site $1$.}
\end{figure}
\end{center}

You can check that, although the hose model is quite simple to specify from the client point of view, it is full of uncertainty for the manager. Indeed, for each source node, for example, the operator ignores how the traffic is shared between the different destination nodes, which constitutes in itself a spatial uncertainty. Furthermore, due to the roughtness of this approach, he does not know how the traffic should evolve under this assumption.
Consequently, we have chosen to model the dynamic evolution of the traffic as a Markov decision process (MDP), which enables us to introduce uncertainty, in our model. 

\bigskip

\small{\textbf{Index of the main notations, used extensively throughout the article.}\\ 
\textit{- $\{X^{(t)}\}_{t\in\mathbb{N}}$- discrete time, discrete state space stochastic process modeling the traffic in the Virtual Private Network $1$.\\
- $X_{ij}^{(t)}$- traffic going from the node i to the node j, at the decision epoch $t$.\\
- $\mathbf{S}$- generic state space.\\
- $\{L^{(t)}\}_{t\in\mathbb{N}}$- traffic on the MPLS network links.\\
- $\mathcal{N}$- set of the sites, or nodes of the MPLS network.\\
- $\mathcal{L}$- set of the links of the MPLS network.\\
- $t_1^{\textrm{out}}$- amount of traffic leaving the site $1$ of the VPN$1$.\\
- $t_2^{\textrm{out}}$- amount of traffic leaving the site $2$ of the VPN$1$.\\
- $t_3^{\textrm{out}}$- amount of traffic leaving the site $3$ of the VPN$1$.\\
- $R(t)$- Routing matrix at the instant t. We note $R$, if the routing is stable, or time invariant.\\
- $\mathbf{S}^X$- discrete state space associated with the Markov Decision Process $\{X^{(t)}\}_{t\in\mathbb{N}}$.\\
- $\mathbf{S}^Y$- discrete state space associated with the Markov Decision Process $\{Y^{(t)}\}_{t\in\mathbb{N}}$.\\
- $\mathbf{S}^Z$- discrete state space associated with the Markov Decision Process $\{Z^{(t)}\}_{t\in\mathbb{N}}$.\\
- $\mathbf{A}$- action space associated with the Markov Decision Process (MDP) $\{X^{(t)}\}_{t\in\mathbb{N}}$.\\
- $F_t(s),\;s\in\mathbf{S}$- vector of strategy associated with the MPD $X^{(t)}$, at the decision epoch t, and for each state $s\in\mathbf{S}\;.$\\
- $f_t(s,a)$- probability for the MDP $\{X^{(t)}\}_{t\in\mathbb{N}}$, to choose the action $a$, in the state $s\in\mathbf{S}$, at the decision epoch t.\\
- $V_t(s)$- value function at the time instant t, in the state $s\in\mathbf{S}$.\\
- $A\;=\;[a_X a_Y a_Z]$- vectors of actions taken at the local level, i.e. on each VPN network.\\
- $a_X\;=\;[a_{X_{12}} a_{X_{21}} a_{X_{31}}]$- actions taken on each link of the VPN$1$.\\ 
- $D\;=\;[d_1 d_2 ... d_{|\mathcal{L}|}]$- actions taken on each link of the whole MPLS network.\\
- $(B^{X}_{ij})^{(t)}$- reserved amount of bandwidth on the directed link $(i,j)$ of the VPN$1$, at time t.\\
- $(B^{Y}_{ij})^{(t)}$- reserved amount of bandwidth on the directed link $(i,j)$ of the VPN$2$, at time t.\\
- $(B^{Z}_{ij})^{(t)}$- reserved amount of bandwidth on the directed link $(i,j)$ of the VPN$3$, at time t.\\
- $(B^{L}_{i})^{(t)}$- reserved amount of bandwidth on the directed link $l_i$ of the global MPLS network, at time t.\\
- $p_{ij}^X$- price associated to the variation of reserved bandwidth on the link (i,j) of the VPN$1$.\\
- $p_{ij}^Y$- price associated to the variation of reserved bandwidth on the link (i,j) of the VPN$2$.\\
- $p_{ij}^Z$- price associated to the variation of reserved bandwidth on the link (i,j) of the VPN$3$.\\
- $\mathbf{Satis_X}$- satisfaction level for the VPN$1$.\\
- $\mathbf{Satis_Y}$- satisfaction level for the VPN$2$.\\
- $\mathbf{Satis_Z}$- satisfaction level for the VPN$3$.\\
- $\mathbf{S}^X\times\mathbf{S}^Y\times\mathbf{S}^Z$- state space of the global process $(X^{(t)},Y^{(t)},Z^{(t)})$.\\
- $\mathbf{A}_{\mathcal{G}}$- action space of the global process $(X^{(t)},Y^{(t)},Z^{(t)})$.\\
- $\mathbf{E}^2$- subset of the global state space $\mathbf{S}^X\times\mathbf{S}^Y\times\mathbf{S}^Z$, where every state violates at least one satisfaction bound.\\
- $\mathbf{E}^1\;=\;\mathbf{S}^X\times\mathbf{S}^Y\times\mathbf{S}^Z\;-\;\mathbf{E}^2\;.$}}

\bigskip 

\section{The representation of Traffic as an MDP}

Let $\{X^{(t)}\}_{t\in\mathbb{N}}$, be the discrete time, discrete state space, stochastic process, representing the traffic in the network. 
The VPN network will be represented by an oriented graph: $G\;=\;(\mathcal{N},\mathcal{L})\;,$ where $\mathcal{N}$ is the set of nodes modeling the sites of the network, and 
$\mathcal{L}$, is the set of directed link of the VPN.\\
Let $X_{ij}^{(t)}$ denotes the traffic going from the node i, to the node j, at the instant t.\\ 
At time period t, the whole traffic is represented by a vector $X^{(t)}$, i.e.:

\begin{equation}
X^{(t)}\;=\;\left( \begin{array}{c}
X^{(t)}_{12}\\
X^{(t)}_{13}\\
\vdots\\
X^{(t)}_{|\mathcal{N}|\;1}\\
\vdots\\
X^{(t)}_{|\mathcal{N}|\;(|\mathcal{N}|-2)}\\
X^{(t)}_{|\mathcal{N}|\;(|\mathcal{N}|-1)}
\end{array} \right)\;.
\end{equation}

The traffic on the link is obtained via the matrix equation:

\begin{equation}
\left( \begin{array}{c}
L^{(t)}_1\\
L^{(t)}_2\\
\vdots\\
L^{(t)}_{|\mathcal{L}|}
\end{array} \right)\;=\;R(t)\;\left( \begin{array}{c}
X^{(t)}_{12}\\
X^{(t)}_{13}\\
\vdots\\
X^{(t)}_{|\mathcal{N}|\;(|\mathcal{N}|-1)}
\end{array} \right)\;,
\end{equation}

where, $R(t)$, models the routing matrix, which can remain constant or change with the time. 
The rather intuitive notation $L^{(t)}_l\;=\;(R(t)\;X^{(t)})_l,\;l\in\mathcal{L}$,  represents the traffic flowing through the link $l$.\\ 
The state space is defined using a simplified version of the hose model. Indeed, the client gives a rather tight upper bound on the traffic going out of each node. As we are supposed to be in the worst case, we should assume that this bound is reached. As an example, in the three-node case, we get a system of relationships:

\begin{equation}
\label{espace_etat}
\left\{ \begin{array}{l}
X_{12}^{(t)}+X_{13}^{(t)}=t_1^{\textrm{out}}\\
X_{21}^{(t)}+X_{23}^{(t)}=t_2^{\textrm{out}}\\
X_{31}^{(t)}+X_{32}^{(t)}=t_3^{\textrm{out}}\\
X_{ij}^{(t)}\geq 0,\;\forall i,j\in\{1,2,3\},\;i\neq j\;.
\end{array} \right.
\end{equation}

We just need to deal with the $3$ components $X_{12}^{(t)},\;X_{21}^{(t)}$ and $X_{31}^{(t)}$, since the others are deduced from the first. The state space is represented geometrically as the union of the three independent segments defined by the system ($\ref{espace_etat}$). Consequently, we will note the continuous state space under the form: $$\mathbf{S}^{\textrm{continuous}}\;=\;\mathbf{S_1}^{\textrm{continuous}}\times\mathbf{S_2}^{\textrm{continuous}}\times\mathbf{S_3}^{\textrm{continuous}}\;.$$
Every element of $\mathbf{S}^{\textrm{continuous}}$ could be represented under a $3$ dimensional vector form: $s\;=\;(s_1,s_2,s_3)$. Where, $s_i$ takes its values in the state space $\mathbf{S_i}^{\textrm{continuous}},\;i=1,2,3\;.$\\
To be more explicit, $\mathbf{S_1}^{\textrm{continuous}}\;=\;\{(X_{12}^{(t)},X_{13}^{(t)})|X_{12}^{(t)}+X_{13}^{(t)}=t_1^{\textrm{out}}\}$, represents the continuous state space associated to the stochastic process $\{X_{12}^{(t)}\}_t$.  $\mathbf{S_2}^{\textrm{continuous}}$  and $\mathbf{S_3}^{\textrm{continuous}}$ define the continuous state spaces associated to the processes $\{X_{21}^{(t)}\}_t$ and $\{X_{31}^{(t)}\}_t$, respectively.
In order to get a discrete state space, the operator should fix \textit{a fiability parameter} $\alpha > 0$, which would characterize the accuracy with which he desires to know the traffic flowing through its links. Then, each of the $3$ segments is discretized using the parameter $\alpha$. The discrete state space resulting, will be logically noted, $\mathbf{S}\;=\;\mathbf{S}_1\times\mathbf{S}_2\times\mathbf{S}_3\;.$
   
\bigskip

\begin{center}
\begin{figure}[h]
\includegraphics[scale=0.5]{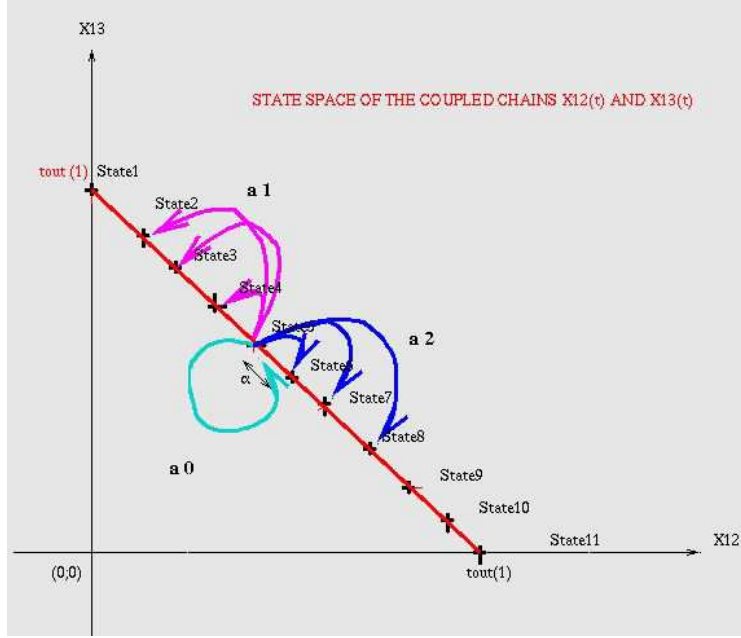}
\caption{Discrete state space and action space.}
\end{figure}
\end{center}

\bigskip

The action space is reduced to $3$ distinct motions, that we will note: $\mathbf{A}\;=\;\{a_0;a_1;a_2\}\;.$ Let us describe the nature of these actions.\\ \\
$\bullet$ If we choose $a_0$, the process stays in the same state.\\
$\bullet$ But, if we choose $a_1$, the traffic increases with an uncertainty on the state transition.
Indeed, we suppose that if the process is in the state $s_i\in\mathbf{S}_j,\;j=1,2,3$ at time t, then it will jump up, on one of the three adjacent states, according to an exponential distribution, decreasing with the distance between these two states. Using the numerotation given in the figure $2$, we simulate a normalized ordered sample of the three transition probabilites $\left[p(s_{i-1}|s_i,a_1),p(s_{i-2}|s_i,a_1),p(s_{i-3}|s_i,a_1)\right]$.
More explicitly,
$$\left\{ \begin{array}{l}
p(s_{i-k}|s_i,a_1)\sim\mathcal{E}(\lambda_1),\;\lambda_1 > 0,\;k=1,2,3,\; \\
p(s_{i-1}|s_i,a_1)\geq p(s_{i-2}|s_i,a_1) \geq p(s_{i-3}|s_i,a_1)\;,
\end{array} \right.$$

under the normalizing constraint: $\displaystyle{\sum_{k=1}^{3}}p(s_{i-k}|s_i,a_1)\;=\;1\;.$ \\
$\mathcal{E}(\lambda_1)$ symbolizes the exponential distribution of parameter $\lambda_1>0\;.$ 
$$X\sim\mathcal{E}(\lambda_1)\;\Leftrightarrow\;f(x;\lambda_1)\;=\;\lambda_1\;\exp^{-\lambda_1\;x}\;\mathbf{1}_{\mathbb{R}^{+}}(x)\;.$$
We notice that if $i=3$, the traffic can only jump on one of the two adjacent states. Consequently, we get the rules:
$$\left\{ \begin{array}{l}
p(s_{i-k}|s_i,a_1)\sim\mathcal{E}(\lambda_1),\;\lambda_1 > 0,\;k=1,2,\; \\
p(s_{i-1}|s_i,a_1)\geq p(s_{i-2}|s_i,a_1)\;,\\
\displaystyle{\sum_{k=1}^{2}}p(s_{i-k}|s_i,a_1)\;=\;1\;.
\end{array} \right.$$
If $i=2$, there is only one possible transition,\\
$p(s_{i-1}|s_i,a_1)\;=\;1\;,$ 
and finally if $i=1$, the traffic has no choice but to stay in the state where it is. This phenomenon results from the finite nature of the state space.\\

$\bullet$ Finally, if we choose $a_2$, the traffic jumps down, on one of the three adjacent states. Formally, we set:
$$\left\{ \begin{array}{l}
p(s_{i+k}|s_i,a_2)\sim\mathcal{E}(\lambda_2),\;\lambda_2 > 0,\;j=1,2,3,\; \\
p(s_{i+1}|s_i,a_2)\geq p(s_{i+2}|s_i,a_2) \geq p(s_{i+3}|s_i,a_2)\;,
\end{array} \right.$$

we still have a normalizing constraint of the form: $\displaystyle{\sum_{k=1}^{3}}p(s_{i+k}|s_i,a_2)\;=\;1\;.$ \\
$\mathcal{E}(\lambda_2)$ symbolizes the exponential distribution of parameter $\lambda_2>0\;.$ 
$$X\sim\mathcal{E}(\lambda_2)\;\Leftrightarrow\;f(x;\lambda_2)\;=\;\lambda_2\;\exp^{-\lambda_2\;x}\;\mathbf{1}_{\mathbb{R}^{+}}(x)\;.$$

If $i\geq (|\mathbf{S}_j|-3)$, we get the same limitations on the transitions as previously mentionned.

\subsection{Iterative game between bandwidth reservation and traffic allocation}

Remind that a strategy specifies for each state $s\in\mathbf{S}$ and each time t, the probability to choose one of the three actions. Under the vector form, we obtain:
$$\forall t\in\mathbb{N},\;\forall s\in\mathbf{S},\;F_t(s)=\left(f_t(s,a_0)\;f_t(s,a_1)\;f_t(s,a_2)\right)^T\;.$$

However, this probability vector is stochastic, and consequently, must satisfy the following constraints of normalisation, and non negativity.

$\left\{ \begin{array}{l}
\displaystyle{\sum_{a\in\mathbf{A}}}f_t(s,a)=1,\;\forall t\in\mathbb{N},\;\forall s\in\mathbf{S},\\
\;f_t(s,a)\geq 0,\forall t\in\mathbb{N},\;\forall s\in\mathbf{S},\;\forall a\in\mathbf{A}\;.
\end{array} \right.$

For each time period t, the strategy is represented by an associated matrix $F_t$.

$$F_t\;=\;(F_t(1)\;F_t(2)\;...\;F_t(N))\;=\;\left( \begin{array}{cccc}
f_t(1,a_0) f_t(2,a_0) \ldots f_t(N,a_0)\\
f_t(1,a_1) f_t(2,a_1) \ldots f_t(N,a_1)\\
f_t(1,a_2) f_t(2,a_2) \ldots f_t(N,a_2)
\end{array} \right)\;.$$

\bigskip

We begin to recall basic definitions, which may be very usefull for a proper understanding of the rest of the article.
\begin{definition}
A strategy is \textbf{stationary}, if it is invariant with respect to the time, i.e.:
$$\forall t\in\mathbb{N},\forall s\in\mathbf{S},\;f_t(s,a)\;=\;f(s,a)\;,\;\forall a\in\mathbf{A}\;,$$
and \textbf{deterministic} or \textbf{pure}, if there exists a unique optimal action for each state, at each instant. Which means that: 
$$\forall s\in\mathbf{S},\;f_t(s,a)\in\{0;1\},\;a\in\mathbf{A}\;.$$\\
\end{definition}
\textbf{At first, we deal with deterministic strategies only. Furthermore, we suppose that the horizon is finite.} \\
We note: $\pi\;=\;(F_0,F_1,...,F_T)$, the sequence of stationary strategies defined on $[0;T]$.\\
To begin with, we consider again the simple model of a $3$ site network. The sites will be numbered 
$$\mathcal{N}\;=\;\{1,2,3\}\;,$$
and are associated with nodes. The directed links are stored in the set 
$$\mathcal{L}\;=\;\{(1,2);(1,3);(2,1);(2,3);(3,1);(3,2)\}\;.$$ Furthermore, we suppose that the routing is stable, i.e. time invariant, and that between each couple of nodes, the only possible path is the directed link joining these two nodes.\\ 

\begin{center}
\begin{figure}[h]
\includegraphics[scale=0.45]{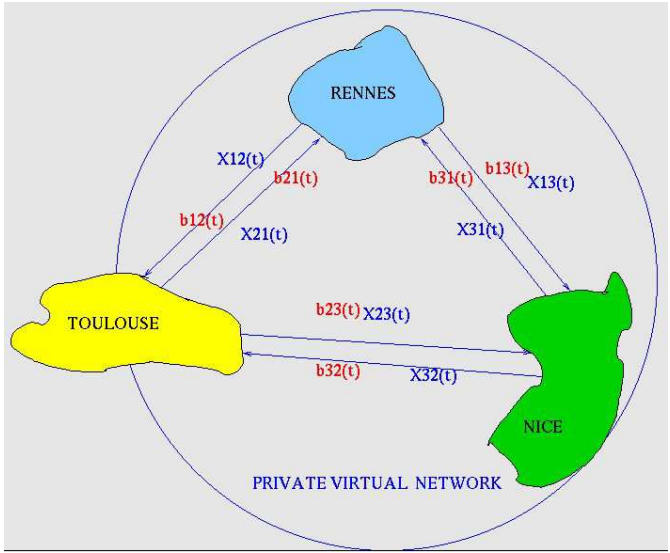}
\caption{A simple model of a $3$ site-VPN.}
\end{figure}
\end{center}

We have chosen to cope with an objective function modeling the delay on the whole network, which is one fundamental parameter in the QoS requirements.
In fact, due to the simple structure of the example, each link is associated with an $M/M/1$ queue, and consequenly the global delay on the whole Virtual Private Network, takes the form:
\begin{equation}
\mathcal{C}_t(X^{(t)})\;=\;\sum_{\{i,j\in\mathcal{N},\;i\neq j\}}\left[\frac{X^{(t)}_{ij}}{B_{ij}^{(t)}-X^{(t)}_{ij}}+p_{ij}(B_{ij}^{(t)}-B_{ij}^{(t-1)})\right]\;.
\end{equation}
Where, $B_{ij}^{(t)}$ is the bandwidth reserved on the directed link $(i,j)$, at time t, by the network manager. 
The second part of this equation stands for a penalty criteria. Indeed, in order to minimize the first part of the equation, the operator should choose to increase infinitely far the reserved amount of bandwidth. Fortunately, the second part, introduces a price $p_{ij}>0$, linked with the variations of the reserved bandwidth. Under this assumption, the manager's interest should be to choose relatively stable values for the amounts of reserved bandwidths.

\bigskip		
	
Our problem takes the formal form:
\begin{equation}
\label{pb_VPN}
\pi^{\star}\;=\;\arg\min_{B=(B^{(0)},B^{(1)},...,B^{(T)})}\arg\max_{\pi=(F_0,F_1,...,F_T)}\left\{\mathbf{E}_{\pi}[\sum_{t=0}^T\beta^t\;\mathcal{C}_t(X^{(t)},F_t)|X_0=s]\;|\;s\in\mathbf{S}\right\}\;.
\end{equation}

\textbf{Remark.} \small{Since the strategies are deterministic, there exists a unique optimal action at each time instant t, and for each state in the state space $\mathbf{S}^X$. Consequently, you guess easily that $F_t(X^{(t)})$, contains the optimal action associated with the random variable $X^{(t)}$, at the instant t. The cost function can then, naturally be interpreted as follows:

\begin{equation}
\mathcal{C}_t(X^{(t)},F_t)\;=\;\sum_{\{i,j\in\mathcal{N},\;i\neq j\}}\left[\frac{X^{(t)}_{ij}+F_t(X^{(t)}_{ij})}{\Phi\left(X^{(t)}_{ij}+F_t(X^{(t)}_{ij})\right)-X^{(t)}_{ij}}+p_{ij}\left(\Phi\left(X^{(t)}_{ij}+F_t(X^{(t)}_{ij})\right)-\Phi\left(X^{(t)}_{ij}\right)\right)\right]\;.
\end{equation}}

The parameter $\beta\in[0;1[$, often called \textit{discount factor}, captures the natural notion that a reward of $1$ unit at a time of $(t+1)$, is worth only $\beta$ of what it was worth at time $t$.\\
In order to simplify the expression of $(\ref{pb_VPN})$, a quite natural idea might be to isolate the sum into two parts. Hence, following our intuition, we write,

$$\pi^{\star}\;=\;\arg\min_{B=(B^{(0)},B^{(1)},...,B^{(T)})}\arg\max_{\pi=(F_0,F_1,...,F_T)}\left\{\mathbf{E}_{\pi}[\mathcal{C}_{0}(X^{(0)},F_0)|X^{(0)}=s]+\mathbf{E}_{\pi}[\sum_{t=1}^T\beta^t\;\mathcal{C}_t(X^{(t)},F_t)|X^{(0)}=s]\right\}\;.$$

Then, it comes easily that,

$$\pi^{\star}\;=\;\arg\min_{B=(B^{(0)},B^{(1)},...,B^{(T)})}\left\{\arg\max_{F_0}[\mathcal{C}_{0}(s,F_0)]+\arg\max_{F_1,F_2,...,F_T}(\mathbf{E}_{\pi}[\sum_{t=1}^T\beta^t\;\mathcal{C}_t(X^{(t)},F_t)|X^{(0)}=s])\right\}\;.$$

If, we repeat once more the same decomposition, we get the expression:

\begin{eqnarray}
\label{toto}
\pi^{\star}\;&=&\;\arg\min_{B=(B^{(0)},B^{(1)},...,B^{(T)})}\{\arg\max_{F_0}[\mathcal{C}_{0}(s,F_0)]+\arg\max_{F_1}\sum_{s'\in\mathbf{S}}[\beta\;\mathcal{C}_1(s',F_1)\;p(s'|s,F_1)]+{}\nonumber\\
{} & &\arg\max_{F_2,F_3,...,F_T}\mathbf{E}_{\pi}[\sum_{t=2}^T\beta^t\;\mathcal{C}_t(X^{(t)},F_t)|X^{(0)}=s]\}\;.
\end{eqnarray}

The equation $(\ref{toto})$ captures the essence of the principle of optimality, which is based on the recursive nature of the equation, and the introduction of the value function $V_t(s),\;s\in\mathbf{S}$.\\
In fact, solving $(\ref{pb_VPN})$ is equivalent to computing a solution to Bellman's optimality equation, which takes the following special setting:

\begin{equation}
V_t(s)\;=\;\min_{B^{(t)}}\;\max_{a\in\mathbf{A}}\left\{\mathcal{C}_{T-t}(s,a)+\beta\;\sum_{s'=1}^{|\mathbf{S}|}p(s'|s,a)\;V_{t-1}(s')\right\},\;\forall s\in\mathbf{S},\;\forall t\in\{0,1,2,...,T\}\;.
\end{equation}

To solve this equation, we proceed by backward induction.

\bigskip

$\star$ To begin with, we suppose that at $t=T$, the reserved bandwidth is fixed. Then \textbf{in each state s}, we have to find the set of actions on each link, which maximizes the equation:
\begin{equation}
\label{fonda}
a_s^{T-1}\;=\;\arg\max_{a\in\mathbf{A}^{|\mathbf{S}|}}\left\{\sum_{i,j}\left[\frac{s_{ij}+a_{ij}}{B_{ij}^{(T)}-(s_{ij}+a_{ij})}+p_{ij}(B_{ij}^{(T)}-B_{ij}^{(T-1)})\right]\right\}\;.
\end{equation}

The traffic being fixed to its new value: $X^{(T)}\;=\;X^{(T-1)}+a^{(T-1)}$, we would like to find the minimal amount of bandwidth to be reserved on each link.
Consequently, we must solve the optimization problem:
\begin{equation}
B^{(T)}\;=\;\arg\min_{B}\left\{\sum_{i,j}\left[\frac{X_{ij}^{(T)}}{B_{ij}-X_{ij}^{(T)}}+p_{ij}(B_{ij}-B_{ij}^{(T-1)})\right]\right\}\;.
\end{equation}

The solution of this continuous optimization problem can be obtained analytically. That's why, we express it as a function of the worst traffic allocation, at time $T$.
\begin{equation}
\label{formal}
B^{(T)}\;=\;\Phi(X^{(T)})\;.
\end{equation}

Finally, substituting $(\ref{formal})$ in the equation $(\ref{fonda})$, we get the simpler expression:
\begin{equation}
a_{x^{(T-1)}}^{T-1}\;=\;\arg\max_{a\in\mathbf{A}^{|\mathbf{S}|}}\left\{\sum_{i,j}\left[\frac{x_{ij}^{(T-1)}+a_{ij}}{\Phi(x_{ij}^{(T-1)}+a_{ij})-(x_{ij}^{(T-1)}+a_{ij})}+p_{ij}(\Phi(x_{ij}^{(T)})-\Phi(x_{ij}^{(T-1)}))\right]\right\}\;,
\end{equation}
 where, $x^{(T-1)}$ is a realization of the traffic process $X^{(T-1)}$ in the state space $\mathbf{S}$.

Now, the value can be easily computed, and we set:
$$V_1(s)\;=\;\mathcal{C}_{T-1}(s,a_s^{T-1}),\;\forall s\in\mathbf{S}\;.$$

$\star$ Then, at the iteration $(T-t),\;t>1$, we proceed exactly the same way. An optimal action, and the associated optimal rewards are known, for the last $(t-1)$ stages. Then, with $t$ stages to go, the only thing we need to do, is to maximize the immediate expected reward and the maximal expected payoff for the remainder of the process with $(t-1)$ stages to go. As a result, we obtain the expression:   

\begin{eqnarray}
a_{x^{(T-t)}}^{T-t}\;&=&\;\arg\max_{a\in\mathbf{A}^{|\mathbf{S}|}}\{\sum_{i,j}\left[\frac{x_{ij}^{(T-t)}+a_{ij}}{\Phi(x_{ij}^{(T-t)}+a_{ij})-
(x_{ij}^{(T-t)}+a_{ij})}+p_{ij}\left(\Phi(x_{ij}^{(T-(t-1))})-\Phi(x_{ij}^{(T-t)})\right)\right]{}\nonumber\\
{} & &+\beta\;\sum_{s_{ij}'=1}^{|\mathbf{S_i}|}p(s_{ij}'|x_{ij}^{(T-t)},a_{ij})\;V_{(t-1)}(s_{ij}')\}\;.
\end{eqnarray}

And, the value takes the form:

\begin{eqnarray}
V_t(s)\;&=&\;\{\sum_{i,j}\left[\frac{x_{ij}^{(T-t)}+a_{ij}^{(T-t)}}{\Phi(x_{ij}^{(T-t)}+a_{ij}^{(T-t)})-
(x_{ij}^{(T-t)}+a_{ij}^{(T-t)})}+p_{ij}\left(\Phi(x_{ij}^{(T-(t-1))})-\Phi(x_{ij}^{(T-t)})\right)\right]{}\nonumber\\
{} & &+\beta\;\sum_{s_{ij}'=1}^{|\mathbf{S_i}|}p(s_{ij}'|x_{ij}^{(T-t)},a_{ij}^{(T-t)})\;V_{(t-1)}(s_{ij}')\}\;.
\end{eqnarray}

We go back with the same idea, until $t=T$.
 
\section{Application to the management of a $3$ site - VPN}

\subsection{Stable and mono-path routing}

To model this simple system, we introduce $3$ MDPs, called respectively $X_1^{(t)}=(X_{12}^{(t)},X_{13}^{(t)})$, $X_2^{(t)}=(X_{21}^{(t)},X_{23}^{(t)})$ and, $X_3^{(t)}=(X_{31}^{(t)},X_{32}^{(t)})$. In fact, these MDPs are not really bi-dimensional, since we only need to take into account the first components. Indeed, the second ones are deduced from the first ones, using the set of equalities ($\ref{espace_etat}$). Since the routing is supposed to be constant, we note that these three MDPs are completely independent of one another. Consequently, our initial problem can be separated into three disjoint sub-problems.

\begin{equation}
\left\{ \begin{array}{l}
 (\pi^{X_1})^{\star}\;=\;\displaystyle{\arg\min_{B^{X_1}}}\;\;\displaystyle{\arg\max_{\pi^{X_1}}}\left\{\mathbf{E}_{\pi^{X_1}}\left[\sum_{t=0}^T\mathcal{C}_t(X_1^{(t)},F_t^{X_1})|X_1^{(0)}=s\right]\;|\;s\in\mathbf{S_1}\right\}\;,\\
 (\pi^{X_2})^{\star}\;=\;\displaystyle{\arg\min_{B^{X_2}}}\;\;\displaystyle{\arg\max_{\pi^{X_2}}}\left\{\mathbf{E}_{\pi^{X_2}}\left[\sum_{t=0}^T\mathcal{C}_t(X_2^{(t)},F_t^{X_2})|X_2^{(0)}=s\right]\;|\;s\in\mathbf{S_2}\right\}\;,\\
 (\pi^{X_3})^{\star}\;=\;\displaystyle{\arg\min_{B^{X_3}}}\;\;\displaystyle{\arg\max_{\pi^{X_3}}}\left\{\mathbf{E}_{\pi_{X_3}}\left[\sum_{t=0}^T\mathcal{C}_t(X_3^{(t)},F_t^{X_3})|X_3^{(0)}=s\right]\;|\;s\in\mathbf{S_3}\right\}\;.\\
 \end{array} \right.
 \end{equation}

To solve these equations, we simply apply Bellman's optimality equation, and backward induction, as explained in the previous section.

\bigskip

%%%%%%%%%%%%%%%%%%%%%%%%%%%%%%%%%%%%%%%%%%%%%%%%%%%%%%%%%%%%%%%%%%%%%%%%%%%%%%%%%%%%%%%%%%%%%%%%%

\begin{center}
\begin{figure}[h]
\includegraphics[scale=0.7]{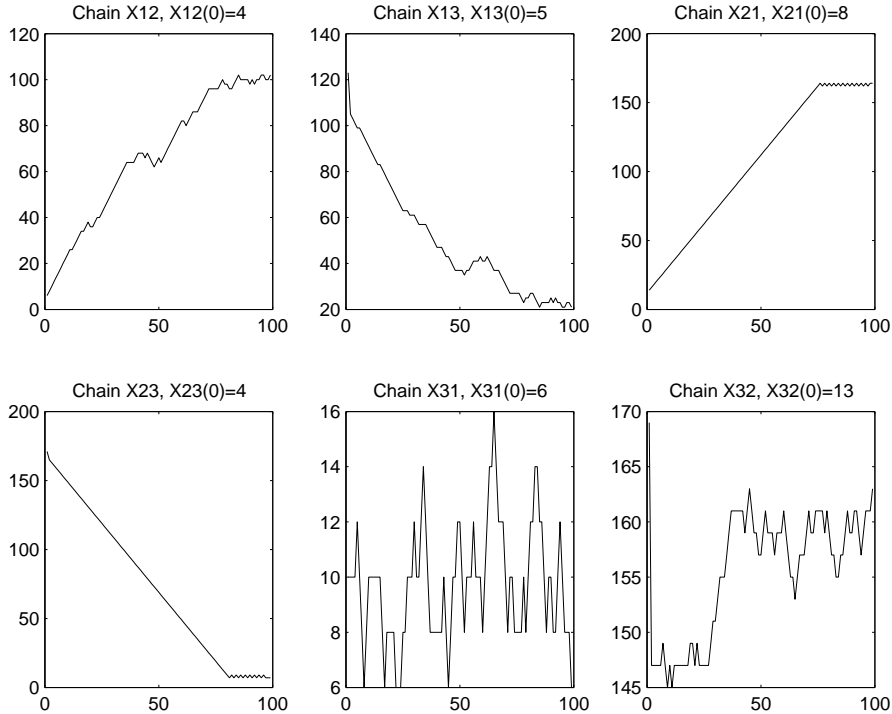}
\caption{Dynamic evolution of the traffic on the six links for $t\leq 100$, $\beta=0.9$.}
\end{figure}
\end{center}

\begin{center}
\begin{figure}[h]
\includegraphics[scale=0.7]{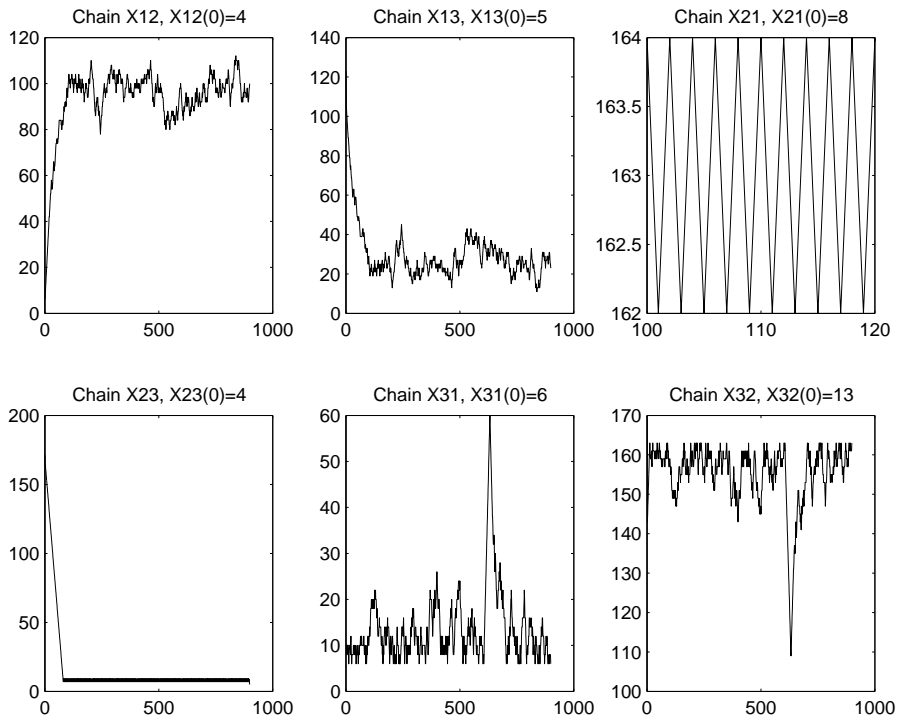}
\caption{A few iterations later ($t\leq 900$, $\beta=0.9$), a stationary behavior appears.}
\end{figure}
\end{center}

%%%%%%%%%%%%%%%%%%%%%%%%%%%%%%%%%%%%%%%%%%%%%%%%%%%%%%%%%%%%%%%%%%%%%%%%%%%%%%%%%%%%%%%%%%%%%%%%%

\bigskip

\begin{center}
\begin{figure}[h]
\includegraphics[scale=0.3]{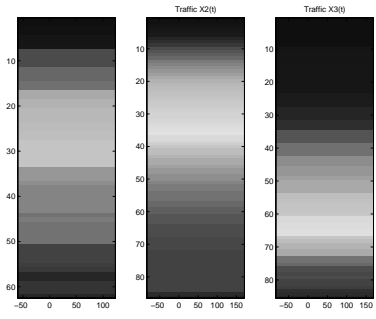}
\caption{Sojourn times for the MDPs $X_{12}^{(t)}$,  $X_{21}^{(t)}$ and $X_{31}^{(t)}$, $t\leq 100$, $\beta=0.9$. The sojourn time associated with a specific couple of state and action $(s,a),\;s\in\mathbf{S},\;a\in\mathbf{A}$, represents the number of times among the decision epochs $\{0,1,...,T\}$, in which we choose the action $a$ in  the state $s$. It enables us to characterize the standard behavior of the system.}
\end{figure}
\end{center}

\begin{center}
\begin{figure}[h]
\includegraphics[scale=0.3]{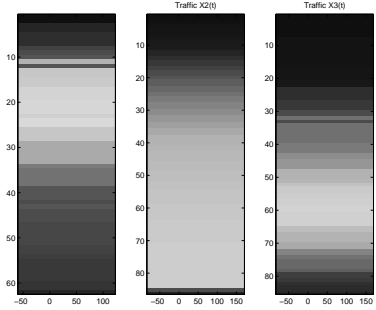}
\caption{Evolution of the sojourn times. The distribution is more homogeneous, since asymptotically the choice of the action in each
 state is time invariant. }
\end{figure}
\end{center}

%%%%%%%%%%%%%%%%%%%%%%%%%%%%%%%%%%%%%%%%%%%%%%%%%%%%%%%%%%%%%%%%%%%%%%%%%%%%%%%%%%%%%%%%%%%%%%%%%

\subsection{Existence of stationary strategies}
The notion of stability is fundamental in the theory of dynamical systems. The usual idea is to prove the convergence of the system towards an equilibrium point. Transposed to the theory of Markov chains, equilibrium points are associated to invariant measures. Finally, in the theory of MDPs, these invariant measures become stationary strategies.\\
We will start by computing the stationary strategies using the elegant approach developed in more details in $[5]$.\\

Let $V$ be an arbitrary vector taking values in the state space $\mathbf{S}$. Using the definition of the optimality equation, we know that in each state $s\in\mathbf{S}$, the value function should satisfy the inequality:

\begin{equation}
\label{opt_in}
V(s)\geq\mathcal{C}(s,a)+\beta\;\sum_{s'=1}^{|\mathbf{S}|}p(s'|s,a)V(s'),\;\forall a\in\mathbf{A},\;\forall s\in\mathbf{S}\;.
\end{equation}

\textbf{Remark.} The function $\mathcal{C}$ is actually independant of the time, since the strategy should be stationary.\\ 

If we multiply the above inequalities by $f(s,a)$, and sum over all the $a\in\mathbf{A}$, we get:

$$V(s)\geq\mathcal{C}(s,F)+\beta\;\sum_{s'=1}^{|\mathbf{S}|}p(s'|s,F)V(s'),\;\forall s\in\mathbf{S}\;,$$

where, $\mathcal{C}(s,F)=\displaystyle{\sum_{a\in\mathbf{A}}}\mathcal{C}(s,a)f(s,a)\;,$ and
$p(s'|s,F)=\displaystyle{\sum_{a\in\mathbf{A}}}p(s'|s,a)f(s,a)$.

In matrix form, the set of inequalities becomes:
$$V\geq\mathcal{C}(F)+\beta\;P(F)V\;,$$

where, the reward function under the strategy F, can be written under the vector form:
$$\mathcal{C}(F)\;=\;\left( \begin{array}{c}
\mathcal{C}(1,F)\\
\mathcal{C}(2,F)\\
\vdots\\
\mathcal{C}(|\mathbf{S}|,F)
\end{array} \right)\;,$$

and the probability transition matrix becomes,

$$P(F)\;=\;\left( \begin{array}{cccc}
p(1|1,F) & p(2|1,F) & \ldots & p(|\mathbf{S}|\;|\;1,F)\\
p(1|2,F) & p(2|2,F) & \ldots & p(|\mathbf{S}|2,F)\\
\vdots & \vdots & \vdots & \vdots\\
p(1\;|\;|\mathbf{S}|,F) & p(2\;|\;|\mathbf{S}|,F) & \ldots & p(|\mathbf{S}|\;|\;|\mathbf{S}|,F)\\
\end{array} \right)\;.$$

Upon substituting the above inequality into itself k times, and taking the limit as $k\rightarrow\infty$, we obtain:

$$V\geq[I-\beta\;P(F)]^{-1}\mathcal{C}(F)\;.$$

But,
$$\sum_{t=0}^{\infty}\beta^t\mathbf{E}_{F}[\mathcal{C}(X^{(t)},F)\;|\;X^{(0)}=s]\\
=\sum_{t=0}^{\infty}\beta^t\;p_t(s'|s,F)\mathcal{C}(s',F)\\
=\sum_{t=0}^{\infty}\beta^t[P^t(F)\;\mathcal{C}(F)\;|\;X^{(0)}=s]\;,$$

since the transition probabilities are time homogeneous.

Going back to our matrix formulation, we get:

$$\left( \begin{array}{c}
\displaystyle{\sum_{t=0}^{\infty}}\beta^t\mathbf{E}_{F}[\mathcal{C}(X^{(t)},F)|X^{(0)}=1]\\
\displaystyle{\sum_{t=0}^{\infty}}\beta^t\mathbf{E}_{F}[\mathcal{C}(X^{(t)},F)|X^{(0)}=2]\\
\vdots\\
\displaystyle{\sum_{t=0}^{\infty}}\beta^t\mathbf{E}_{F}[\mathcal{C}(X^{(t)},F)|X^{(0)}=|\mathbf{S}|]
\end{array} \right)\;=\;[I-\beta\;P(F)]^{-1}\mathcal{C}(F)\;.$$

We see that an arbitrary vector $V$ satsfying $(\ref{opt_in})$, is an upper bound on the discounted value vector due to any stationary strategy F. Consequently, we should naturally think that the discounted value vector might be the optimal solution of the linear programm:
$$\left\{ \begin{array}{l}
\min\;\displaystyle{\sum_{s=1}^{|\mathbf{S}|}}\gamma(s)V(s)\\
V(s)\geq\mathcal{C}(s,a)+\beta\;\displaystyle{\sum_{s'=1}^{|\mathbf{S}|}}p(s'|s,a)V(s'),\;a\in\mathbf{A},\;s\in\mathbf{S}
\end{array} \right.$$

with $\gamma(s)$ ($\gamma(s)>0$ and $\displaystyle{\sum_{s\in\mathbf{S}}}\gamma(s)=1$), being the probability that the process begins in state $s\in\mathbf{S}$.\\ 
By duality, we get:

\begin{equation}
\left\{ \begin{array}{l}
\label{LP}
\max\displaystyle{\sum_{s=1}^{|\mathbf{S}|}}\sum_{a=1}^{3}\mathcal{C}(s,a)\;x_{sa}\\
\displaystyle{\sum_{s=1}^{|\mathbf{S}|}}\sum_{a=1}^{3}[\delta(s,s')-\beta p(s'|s,a)]x_{sa}\;=\;\gamma(s),\;s'\in\mathbf{S}\\
x_{sa}\geq 0,\;a\in\mathbf{A},\;s\in\mathbf{S}\;.
\end{array} \right.
\end{equation}

$x_{sa},\;s\in\mathbf{S},\;a\in\mathbf{A}$, which is the dual variable, can be heuristically interpreted as the long-run fraction of decision epochs at which the system is in the state s, and the action a is made.\\

\textbf{Remark.} It can be shown (see $[5]$), that if the variables $x_{sa}$ are obtained using the simplex algorithm, as the solution of the above linear program, then the associated stationary strategy is indeed a deterministic stationary strategy.
We note, $x_s\;=\;\displaystyle{\sum_{a\in\mathbf{A}}}x_{sa}\;,$ for each $s\in\mathbf{S}$. In our context, we get the following system of equalities:

\begin{equation}
\label{simplex}
\left\{ \begin{array}{l}
x_{1 a_0}+x_{1 a_1}+x_{1 a_2}\;=\;x_1\\
x_{2 a_0}+x_{2 a_1}+x_{2 a_2}\;=\;x_2\\
\vdots\\
x_{|\mathbf{S}| a_0}+x_{|\mathbf{S}| a_1}+x_{|\mathbf{S}| a_2}\;=\;x_{|\mathbf{S}|}\\
x_{sa}\geq 0,\;\forall s\in\mathbf{S},\;\forall a\in\mathbf{A}\;.
\end{array} \right.
\end{equation}

However, since $x^0=\left(x_{1 a_O}^0 x_{1 a_1}^0 x_{1 a_2}^0|x_{2 a_0}^0 x_{2 a_1}^0 x_{2 a_2}^0|\ldots|x_{|\mathbf{S}| a_0}^0 x_{|\mathbf{S}| a_1}^0 x_{|\mathbf{S}| a_2}^0\right)$, is an optimal basic feasible solution for the simplex algorithm, it is necessary an extreme point of the space defined by the system $(\ref{simplex})$. Using the definition of an extreme point, for all $s\in\mathbf{S}$, there exists a unique $i\in\{0;1;2\}$ such that $x_{s a_i}^0=x_s^0$, and $x_{s a_j}^0=0$, for $j\in\{0;1;2\},j\neq i$. Consequently, the stationary control $F^0$ constructed from $x^0$, by setting:
$$f^0(s,a)\;=\;\frac{x_{sa}^0}{\displaystyle{\sum_{a\in\mathbf{A}}}x_{sa}^0}\;,\;\forall a\in\mathbf{A},\;\forall s\in\mathbf{S},$$ is deterministic.

\bigskip

We have proved the existence of a stationary strategy for our problem. But, are we sure that the set of feasible strategies determined on $[0;T]$ with the help of Bellman's optimality equation converges asymptotically to these values?
The concept of ergodicity is well-known in the field of system engineering, and more specifically in queuing theory. The problem with such systems, is to define operational parameters which should optimize the performance of our system. Those parameters are usually deduced from the study of the stationary behavior of the global system. 
The idea should be to study the dynamic evolution of a given trajectory of the system. But, do all these specific realizations adopt the same asymptotic behavior?
Is there a law linking operational and stochastic performance parameters?\\ 
Basically, \textit{a system is said to be ergodic}, if all the specific realizations of the dynamic evolution of the system are asymptotically and statistically the same. In fact, ergodicity is synonimous with equality between spatial and temporal means. As a result, in such a framework, the operational parameters are equal to the stochastic performance parameters. 
Translated to the MDP context, the property of ergodicity is defined conditionally on the choice of an action.   
The theoretic definition below introduces the notion of ergodicity from a measure theoretic point of view.

\begin{defn}
Conditionally to the choice of an action $a\in\mathbf{A}$, the Markov chain\\ 
$(\mathbf{S},\{p(s'|s,a)\}_{s,s'\in\mathbf{S}})$ is ergodic if:
\begin{equation}
|p[X^{(t+1)}\in\mathbf{B}|X^{(0)}=s]-\mu^{\star}(\mathbf{B})|\rightarrow 0,\;\forall s\in\mathbf{S},\;\forall \mathbf{B}\in\mathcal{B}(\mathbf{S})\;,
\end{equation}
\end{defn}

where, $\mu^{\star}(\mathbf{B})\;=\;\int\mu(dx)\;\mathbf{P}[X^{(t+1)}\in\mathbf{B}\;|\;X^{(0)}=x]\;.$ 
Besides, a well known result states that in the case of an ergodic MDP, the set of strategies $\{f_t(s,a)|s\in\mathbf{S},\;a\in\mathbf{A}\}_t$, converges to a set of stationary strategies $\{f(s,a)|s\in\mathbf{S},\;a\in\mathbf{A}\}$, i.e. strategies which are time invariant.\\
In practice, we would rather use the fudamental result evoked before, which states that to prove the ergodicity of a Markov chain, it suffices to establish the equality between temporal and spatial order means. 
Since conditionally to the choice of an action $(\mathbf{S},\{p(s'|s,a)\}_{s,s'\in\mathbf{S}})$ is a Markov chain, we are able to transpose this result to the theory of MDP. 
Then, for every $a\in\mathbf{A}$, we have to verify that:
\begin{equation}
\forall k\in\mathbb{N},\;\displaystyle{\lim_{t\rightarrow\infty}}\sum_{s\in\mathbf{S}}s^k\;f_t(s,a)\;=\;\displaystyle{\lim_{T\rightarrow\infty}}\sum_{s\in\mathbf{S}}s^k\frac{\hat{x_{sa}}}{T}\;,\;\textit{almost everywhere}\;.
\end{equation}

 The second part of the equality is obtained using simulation. Indeed, using Bellman's optimality principle, and for $T$ large enough, we are given a set of optimal deterministic strategies on $[0;T]$: $\pi\;=\;\left(F_0,F_1,...,F_T\right)\;.$ In order to build sample trajectories, we just have to choose an arbitrary initial state, or even better, an initial distribution on the state space. Then, we should find the optimal action associated to the state, at the decision epoch t. As a result of this action, we are driven in a new state, and we repeat the process until $t=T$.\\
To compute the first part of the equation, since the first term is positive, we can interchange the sum and the limit.
And, to determine $f(s,a)$, we use linear programming. 
Indeed, the linear program $(\ref{LP})$, introduced previously, gives us the values of the optimal parameters $\{x_{sa}\;|\;s\in\mathbf{S},\;a\in\mathbf{A}\}$, via the simplex algorithm. 
The strategies obtained using the normalizing ratio,
$$f(s,a)\;=\;\frac{x_{sa}}{\displaystyle{\sum_{a\in\mathbf{A}}}x_{sa}}\;,\;\forall s\in\mathbf{S},\;\forall a\in\mathbf{A},$$
are stationary and deterministic. 
Since the empirical and the statistical means coincide asymptotically, we deduce that the sequence of optimal strategies determined by dynamic programming, converges to the stationary strategy obtained using linear programming.

\bigskip

\begin{center}
\begin{tabular}{c}
\includegraphics[scale=0.54]{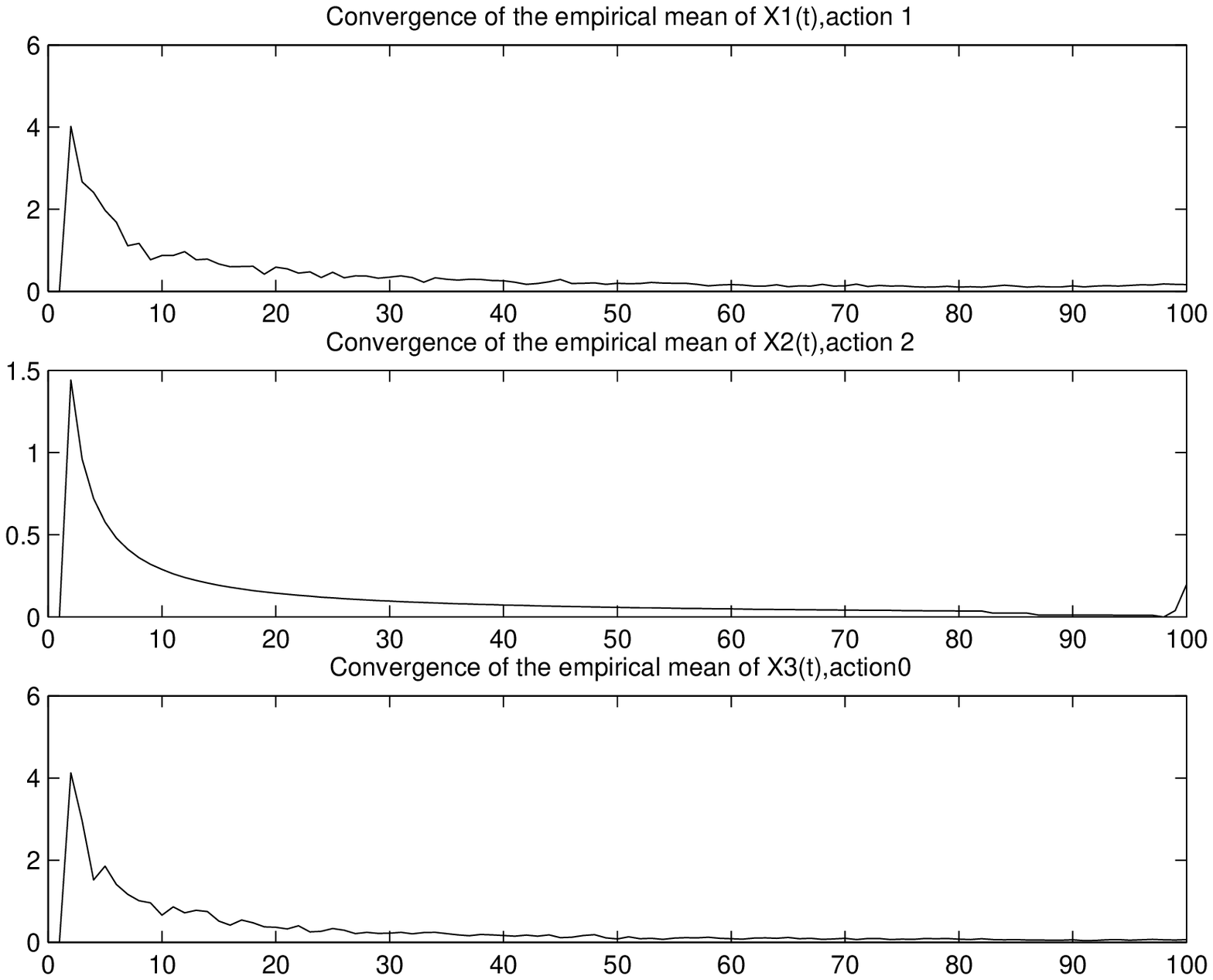}\\
\small{Convergence of the empirical means of the MDPs conditionally to the choice of the actions, $t\leq 100$.}\\
\includegraphics[scale=0.54]{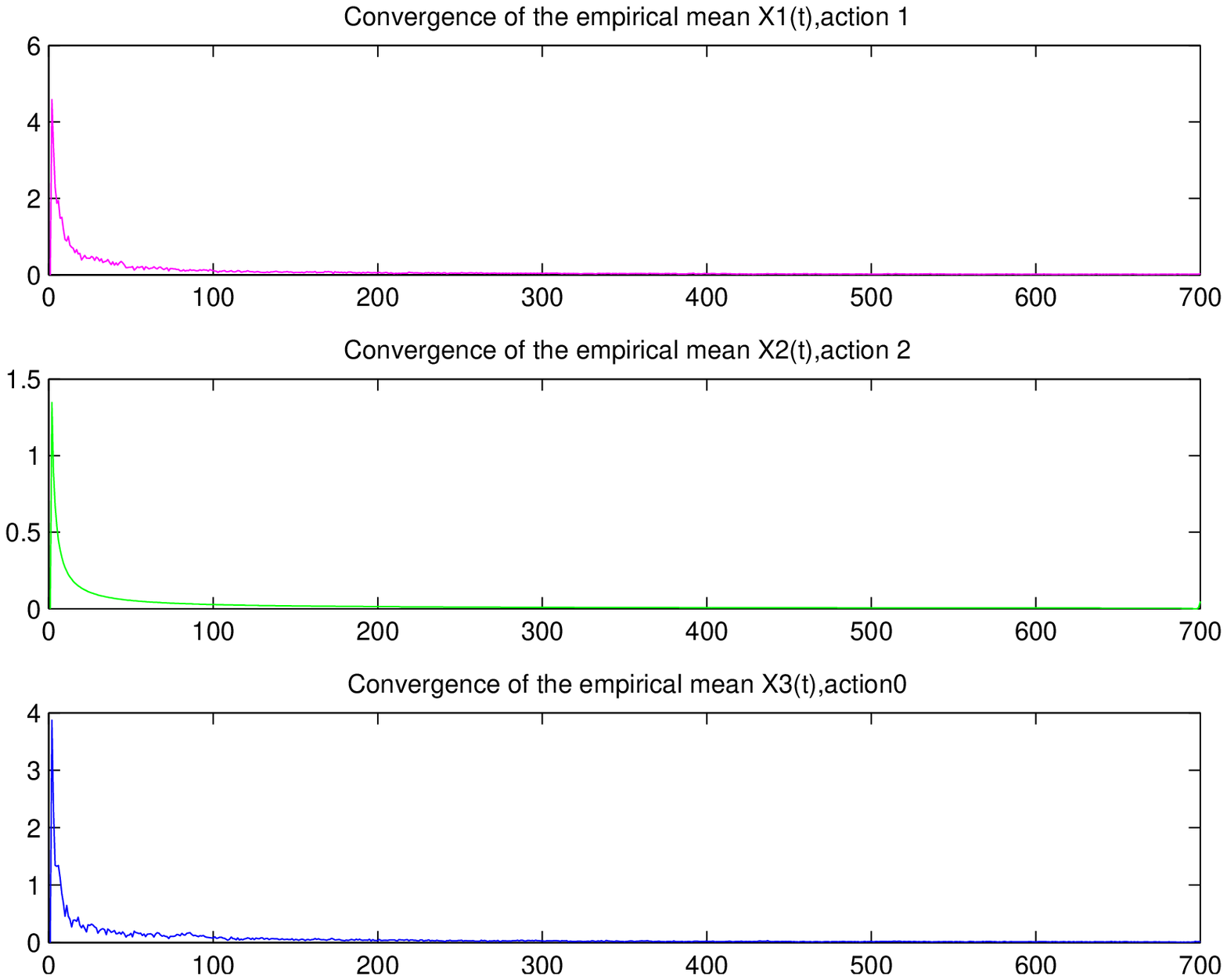}\\
\small{Convergence of the empirical means of the MDPs conditionally to the choice of the actions, $t\leq 300$.}\\
\end{tabular}
\end{center}

\bigskip

\textbf{Remark.} What happen's to this model, if we suppose that the routing is multi-path and changing? We notice that the three MDPs are not independent anymore. 
Consequently, the state space and the action space are made of all the combinations of $3$ elements taken from the initial state space $\mathbf{S}^X$, and the initial action space $\mathbf{A}$, respectively. 
Let $\mathcal{L}$, be the set of links of the network, and $R(t)$, the routing matrix at time t.
The optimality equation remains unchanged on the form, but the cost function is more complicated.

\begin{equation}
V_t(s)\;=\;\min_{B^{(t)}}\;\max_{a\in\mathbf{A}}\left\{\mathcal{C}_{T-t}(s,a)+\sum_{s'=1}^{|\mathbf{S}|}p(s'|s,a)\;V_{t+1}(s')\right\},\;\forall s\in\mathbf{S},\;\forall t\in\{0,1,2,...,T\}\;,
\end{equation}

where, $$C_t(X^{(t)})\;=\;\sum_{l=1}^{|\mathcal{L}|}\left[\frac{(R(t)\;X^{(t)})_l}{B_l^{(t)}-(R(t)\;X^{(t)})_l}+p_l\left((R(t)\;X^{(t)})_l-(R(t-1)\;X^{(t-1)})_l\right)\right]\;,$$

$p_l > 0$, is the price associated with the link $l$.\\
Besides, conditionally to the choice of a three-dimensional action, we make the assumption that the transition probabilities are independant of one another, i.e.: 
$$p\left[(s_1',s_2',s_3')\;|\;(s_1,s_2,s_3),(a_1,a_2,a_3)\right]\;=\;p\left[s_1'\;|\;(s_1,a_1)\right]\;p\left[s_2'\;|\;(s_2,a_2)\right]\;p\left[s_3'\;|\;(s_3,a_3)\right]\;,\;\forall (s_1,s_2,s_3),\;(s_1',s_2',s_3')\in\mathbf{S}.$$

\begin{center}
\begin{tabular}{c}
\includegraphics[scale=0.54]{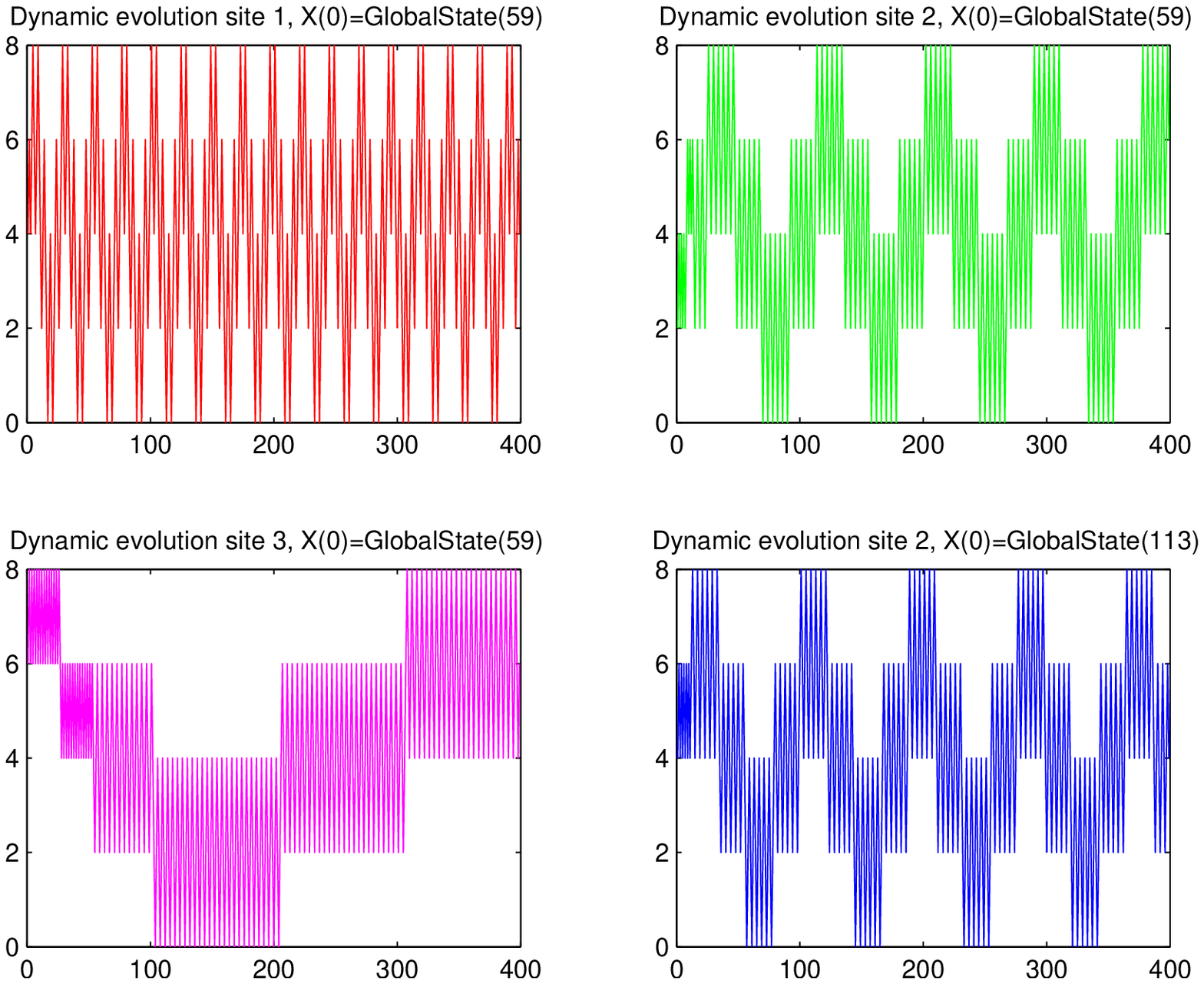}\\
\small{Dynamic evolution of the traffic on the VPN for a changing routing, with exponentialy distributed weights.}\\
\includegraphics[scale=0.54]{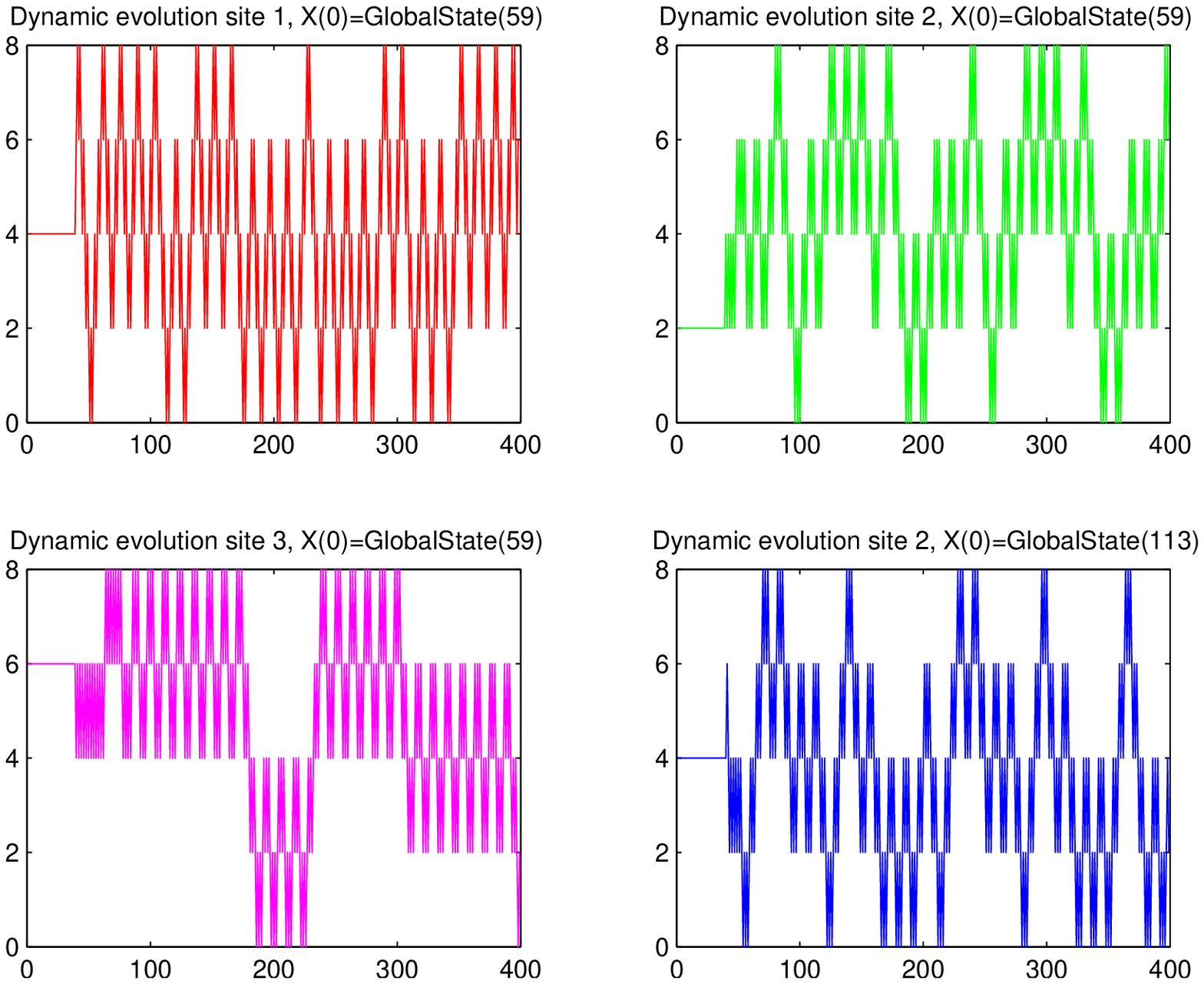}\\
\small{Dynamic evolution of the traffic on the VPN for a changing routing, with normaly distributed weights.}\\
\end{tabular}
\end{center}

\section{The need for a centralized control in a MPLS-network}

We have introduced a way to compute the worst dynamic evolution of the traffic on a $3$ site-VPN with or without changing routing. Going deeper in the reflexion, we can question ourselves about the possibility to manage a network of at least $3$ distinct VPNs.

\begin{center}
\begin{figure}[h]
\includegraphics[scale=0.5]{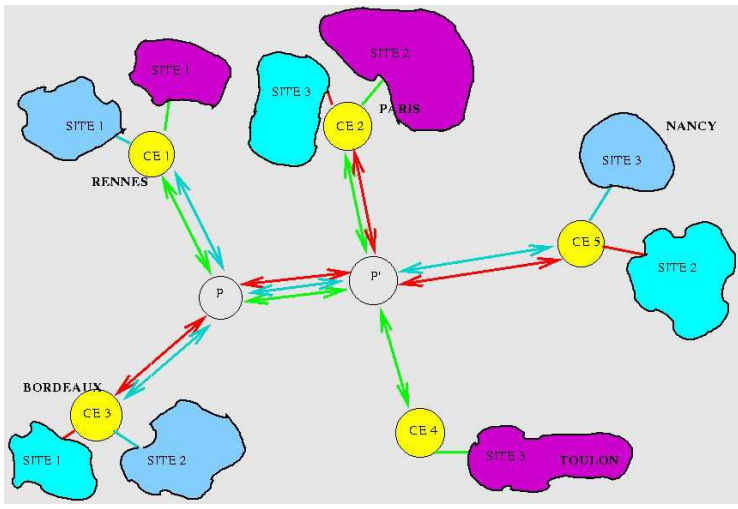}
\caption{A MPLS Network of $3$ VPNs.}
\end{figure}
\end{center}

Each VPN expects that the manager would satisfy the level of QoS it has chosen in the SLA. 
Consequently, the operator should be able to manage $3$ different VPNs, sharing the same infrastructure, each having specific requirements on the quality of service. Furthermore, the VPNs are not aware of one another presence in the network. In fact, at the local level, each VPN behaves completly selfishly, insofar as it tries to optimize its own criteria, using the shared bandwidth, without any knowledge of the needs of the others.\\ \\
Let $\mathcal{L}\;=\;\{1,2,..,12\}$, be the set of the links on the global network.\\ \\
$\mathbf{S}^l\;=\;\{1,2,...,|\mathbf{S}^l|\},\;l\in\mathcal{L}$, contains the state space associated with the link $l$, in the MPLS network.\\ \\
Actually, we define a bound for each VPN, which models the admissible level of delay that the client is able to bear.
We note these levels:
 $\mathbf{Satis}_{X},\;\mathbf{Satis}_{Y}$ and, $\mathbf{Satis}_{Z}$, respectively.\\ \\
The stochastic process associated to the VPN$1$, will be noted:
$X^{(t)}=(X_{12}^{(t)},X_{21}^{(t)},X_{31}^{(t)})$.\\
Recall that $X_{ij}^{(t)}$, stands for the traffic going out of the site i towards the site j of the VPN$1$, at the instant t.
Similarly, the process associated with the VPN$2$ is denoted:
$Y^{(t)}=(Y_{12}^{(t)},Y_{21}^{(t)},Y_{31}^{(t)})$.\\
$Y_{ij}^{(t)}$ represents the traffic flowing from the site i, to the site j, on the VPN$2$, at the instant t.
Finally, $Z^{(t)}=(Z_{12}^{(t)},Z_{21}^{(t)},Z_{31}^{(t)})$, will represent the traffic on the third VPN.\\ 
At the local level, the game is still played the same. At each time, the operator makes a bandwidth reservation on the link of the VPN, the traffic chooses the worst associated allocation on the VPN's links. At the global level, the decisions are centralized. Indeed, the actions are chosen directly on the links of the global MPLS network. 
As in  the local approach, there are three distinct available actions for each link: 
$$\mathbf{D}\;=\;\{d^0;d^1;d^2\}\;.$$
$\bullet$ $d^0$, means that the traffic on the link remains unchanged,\\ \\
$\bullet$ $d^1$, means that the traffic on the link increases with some uncertainty on the next state it will enter.\\ 
$\left\{ \begin{array}{l}
p(s_{i-j}|s_i,d^1)\sim\mathcal{E}(\nu_1),\;\nu_1>0,\;j=1,2,3,\\
p(s_{i-1}|s_i,d^1)\geq p(s_{i-2}|s_i,d^1)\geq p(s_{i-3}|s_i,d^1)\;,\;s_i\in\mathbf{S}^l\;,
\end{array}\right.$ \\ 
under the normalizing constraint: $\displaystyle{\sum_{k=1}^{3}}p(s_{i-k}|s_i,d^1)=1$, and with the the usual limitations, due to the finite cardinality of the state space.\\ \\
$\bullet$ $d^2$, means that the traffic on the link decreases, and the laws are of the same type as previously explained:\\
$\left\{ \begin{array}{l}
p(s_{i+j}|s_i,d^2)\sim\mathcal{E}(\nu_2),\;\nu_2>0,\;j=1,2,3,\\
p(s_{i+1}|s_i,d^2)\geq p(s_{i+2}|s_i,d^2)\geq p(s_{i+3}|s_i,d^2),\;s_i\in\mathbf{S}^l\;,
\end{array}\right.$ \\ \\
with: $\displaystyle{\sum_{k=1}^{3}}p(s_{i+k}|s_i,d^2)=1\;.$\\
The idea now, is to define a rule, which would give the optimal decision epochs at which the decisions should be taken centrally, while the control should be chosen locally during the rest of the time. Consequently, we choose the following rule :\\

\textbf{Rule.} \textit{If one of the bound is not satisfied, then the decisions are taken at the global level, i.e. on the link of the \textbf{whole} network, until all the bounds become satisfied.} \\
 
\begin{center}
\begin{figure}[h]
\includegraphics[scale=0.5]{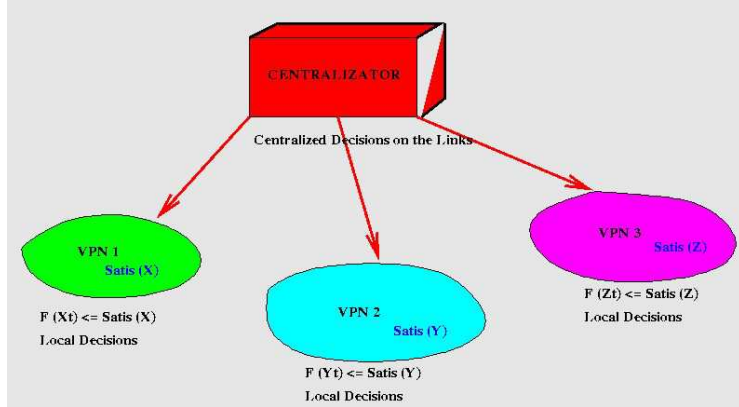}
\caption{A centrally managed network.}
\end{figure}
\end{center}

\subsection{A hierarchical finite horizon MDP approach}

In this section, we make the assumption that the horizon is finite. \\
We start by defining a global process on the MPLS network: $(X^{(t)},Y^{(t)},Z^{(t)})$.\\
- $X^{(t)}\in\mathbf{S}^X$, represents the amount of traffic on the VPN$1$,\\
- $Y^{(t)}\in\mathbf{S}^Y$, is the volume of traffic flowing through the VPN$2$,\\
- $Z^{(t)}\in\mathbf{S}^Z$, is the amount of traffic on the VPN$3$.\\
- $L^{(t)}\;=\;R\;\left( \begin{array}{c}
X^{(t)}\\
\hline\\
Y^{(t)}\\
\hline\\
Z^{(t)}
\end{array} \right)\;,$ represents the amounts of traffic on the links of the MPLS network. $R$, is the routing matrix, which will be supposed to stay stable, i.e. there is no change in the routing. 

\bigskip

The actions at the local level, will be denoted: 
$$A(s(t))\;=\;\left[a_X(x(t)),a_Y(y(t)),a_Z(z(t))\right],\;\forall (x(t),y(t),z(t))\in\mathbf{S}^X\times\mathbf{S}^Y\times\mathbf{S}^Z,\;\forall\;t\in\{0,1,...,T\}\;.$$ 
Indeed, each components of this vector is associated with the action that should be taken on each site of each VPN, at the instant t, in each possible state.
For example, $a_X(x(t))\;=\;\left[a_{X_{12}}(x_{12}(t)),a_{X_{21}}(x_{21}(t)),a_{X_{31}}(x_{31}(t))\right]$ contains the actions to be taken on the links $(1,2)$, $(2,1)$, and $(3,1)$ respectively, of the VPN$1$.\\
At the global level, the actions are centralized, and noted: 
$$D(l(t))\;=\;\left[d_1(l_1(t)),d_2(l_2(t)),...,d_{|\mathcal{L}|}(l_{|\mathcal{L}|}(t))\right],\;\forall\;l_i(t)\in\mathbf{S}^i,\;\forall t\in\{0,1,...,T\}\;.$$
Each element $d_i(l_i(t))$, stands for the action to be taken on the link i, of the MPLS network, at time t, provided the link is in the specific state $l_i(t)$.

\bigskip

At the local level, we have to determine one optimal sequence of strategies per VPN. Formally, the problem can be written under the form:
\begin{equation}
\left\{ \begin{array}{l}
 (\pi^X)^{\star}\;=\;\displaystyle{\arg\min_{B^X}}\;\;\displaystyle{\arg\max_{\pi^X}}\left\{\mathbf{E}_{\pi^X}\left[\sum_{t=0}^T\mathcal{C}_t^X(X^{(t)},F_t^X)|X^{(0)}=s_1\right]\;|\;s_1\in\mathbf{S^X}\right\}\;,\\
 (\pi^Y)^{\star}\;=\;\displaystyle{\arg \min_{B^Y}}\;\;\displaystyle{\arg\max_{\pi^Y}}\left\{\mathbf{E}_{\pi^Y}\left[\sum_{t=0}^T\mathcal{C}_t^Y(Y^{(t)},F_t^Y)|Y^{(0)}=s_2\right]\;|\;s_2\in\mathbf{S^Y}\right\}\;,\\
 (\pi^Z)^{\star}\;=\;\displaystyle{\arg \min_{B^Z}}\;\;\displaystyle{\arg\max_{\pi^Z}}\left\{\mathbf{E}_{\pi^Z}\left[\sum_{t=0}^T\mathcal{C}_t^Z(Z^{(t)},F_t^Z)|Z^{(0)}=s_3\right]\;|\;s_3\in\mathbf{S^Z}\right\}\;,\\
 \end{array} \right.
 \end{equation}

 On the links, the structure of the decision problem is the same, except that the decisions are chosen on each link separately.
 \begin{equation}
 (\pi^L)^{\star}\;=\;\displaystyle{\arg\min_{B^L}}\;\;\displaystyle{\arg\max_{\pi^L}}\left\{\mathbf{E}_{\pi^L}\left[\sum_{t=0}^T\mathcal{C}_t^L(L^{(t)},F_t^L)|L^{(0)}=l\right]\;|\;l\in\mathbf{S^1}\times\mathbf{S^2}\times...
 \times\mathbf{S^{|\mathcal{L}|}}\right\}\;.
 \end{equation}
 
 Recall that the reward functions, defined for each of our $3$ site-VPN, are of the form:
 \begin{equation}
\left\{ \begin{array}{l}
\mathcal{C}^X_t(X^{(t)})\;=\;\displaystyle{\sum_{\{i,j\in\{1,2,3\},\;i\neq j\}}}\left[\frac{X^{(t)}_{ij}}{(B_{ij}^X)^{(t)}-X^{(t)}_{ij}}+p^X_{ij}\left((B_{ij}^X)^{(t)}-(B_{ij}^X)^{(t-1)}\right)\right]\;,\\
\mathcal{C}^Y_t(Y^{(t)})\;=\;\displaystyle{\sum_{\{i,j\in\{1,2,3\},\;i\neq j\}}}\left[\frac{Y^{(t)}_{ij}}{(B_{ij}^Y)^{(t)}-Y^{(t)}_{ij}}+p^Y_{ij}\left((B_{ij}^Y)^{(t)}-(B_{ij}^Y)^{(t-1)}\right)\right]\;,\\
\mathcal{C}^Z_t(Z^{(t)})\;=\;\displaystyle{\sum_{\{i,j\in\{1,2,3\},\;i\neq j\}}}\left[\frac{Z^{(t)}_{ij}}{(B_{ij}^Z)^{(t)}-Z^{(t)}_{ij}}+p^Z_{ij}\left((B_{ij}^Z)^{(t)}-(B_{ij}^Z)^{(t-1)}\right)\right]\;.\\
\end{array} \right.
 \end{equation}
 
According to the same idea, the reward on the links is of the type:
 \begin{equation}
 \mathcal{C}^L_t(L^{(t)})\;=\;\sum_{\{i\in\{1,2,...,|\mathcal{L}|\}\}}\left[\frac{L^{(t)}_{i}}{(B_{i}^L)^{(t)}-L^{(t)}_{i}}+p^L_{i}\left((B_{i}^L)^{(t)}-(B_{i}^L)^{(t-1)}\right)\right]\;.
 \end{equation}
 We then put Bellman's principle of optimality in application, in order to get the sequences of optimal deterministic strategies: $(\pi^X)^{\star},\;(\pi^Y)^{\star},\;(\pi^Z)^{\star}$, and $(\pi^L)^{\star}$.\\ 
 Suppose now, that we have computed these optimal strategies. The next problem we have to face, is how we could build the optimal trajectories of the worst traffic process, for each VPN.\\
 We begin to choose an initial state for each trajectories. \\
 $\star$ $(t=0)$.\\ 
 -On the VPN$1$, we will note $x^{(0)}\;=\;(x^{(0)}_{12},x^{(0)}_{21},x^{(0)}_{31})$, the chosen initial state,\\ 
 -on the VPN$2$, we will choose: $y^{(0)}\;=\;(y^{(0)}_{12},y^{(0)}_{21},y^{(0)}_{31})$,\\
 -and finally, on the VPN$3$, we let: $z^{(0)}\;=\;(z^{(0)}_{12},z^{(0)}_{21},z^{(0)}_{31})$.
 
 \bigskip
 
$\star$ $(t=1)$.\\
If $\mathcal{C}(x^{(0)})\leq \textrm{Satis}_{X}$, and, $\mathcal{C}(y^{(0)})\leq \textrm{Satis}_{Y}$, and, $\mathcal{C}(z^{(0)})\leq \textrm{Satis}_{Z}$, then, we choose the associated optimal actions, and get:
 
 $$\left\{ \begin{array}{l}
 x^{(1)}\;=\;x^{(0)}+a_X(x^{(0)})\;,\\
 y^{(1)}\;=\;y^{(0)}+a_Y(y^{(0)})\;,\\
 z^{(1)}\;=\;z^{(0)}+a_Z(z^{(0)})\;.
\end{array} \right.$$
 
$\star$  At the t$^{th}$ iteration, we check whether or not, $\mathcal{C}(x^{(t)})\leq \textrm{Satis}_{X}$, and, $\mathcal{C}(y^{(t)})\leq \textrm{Satis}_{Y}$, and, $\mathcal{C}(z^{(t)})\leq \textrm{Satis}_{Z}$.\\
$\star\star$ If it is the case, we follow exactly the same way, and obtain:
$$\left\{ \begin{array}{l}
 x^{(t+1)}\;=\;x^{(t)}+a_X(x^{(t)})\;,\\
 y^{(t+1)}\;=\;y^{(t)}+a_Y(y^{(t)})\;,\\
 z^{(t+1)}\;=\;z^{(t)}+a_Z(z^{(t)})\;.
\end{array} \right.$$

$\star\star$ However, if the levels are overwhelmed, then, the decisions are centralized. We start by computing the associated amount of traffic on each link of the MPLS network.
In matrix form, we get:

\begin{equation}
\left( \begin{array}{c}
l_1^{(t)}\\
l_2^{(t)}\\
\vdots\\
l_{|\mathcal{L}|}^{(t)}
\end{array} \right)\;=\;R\;\left( \begin{array}{c}
x^{(t)}\\
\hline\\
y^{(t)}\\
\hline\\
z^{(t)}
\end{array} \right)\;.
\end{equation}

As we actually know in which state the MPLS network globally lies, we choose the optimal action associated. This action tells us the worst way the traffic behaves on each link of the MPLS network. 

$$\left( \begin{array}{c}
l_1^{(t+1)}\\
l_2^{(t+1)}\\
\vdots\\
l_{|\mathcal{L}|}^{(t+1)}
\end{array} \right)\;=\;\left( \begin{array}{c}
l_1^{(t)}+d_1(l_1^{(t)})\\
l_2^{(t)}+d_2(l_2^{(t)})\\
\vdots\\
l_{|\mathcal{L}|}^{(t)}+d_{|\mathcal{L}|}(l_{|\mathcal{L}|}^{(t)})
\end{array} \right)\;.$$

$\star$ At $(t+1)$, we have to check whether or not the levels are satisfied. But, we only know the global amounts of traffic on each link of the MPLS network. In fact, we need to determine the amounts of traffic flowing through each oriented couple of nodes, on each VPN. The traffic being model as a global matrix for each VPN, we have to cope with the matrix equation:

\begin{equation}
\label{MT}
\left( \begin{array}{c}
l_1^{(t+1)}\\
l_2^{(t+1)}\\
\vdots\\
l_{|\mathcal{L}|}^{(t+1)}
\end{array} \right)\;=\;R\;\left( \begin{array}{c}
x^{(t+1)}\\
\hline\\
y^{(t+1)}\\
\hline\\
z^{(t+1)}
\end{array} \right)\;.
\end{equation}

Unfortunately, the problem is severly undertermined, in most applications. 

\subsection{How to jump from a global level to local levels?}
Various statistical techniques of estimation can be employed to solve such problems. In this paper, we have chosen to use an original approach, at least in this field, based on the Cross-Entropy method ($[16]$). Indeed, this technique seems to be well-adapted to solve problems of changing routing, and consequently, it could be envisaged to be used in extensions of our approach.

\subsubsection{A brief introduction to the Cross-Entropy (CE) method}
The CE method ($[16]$), is a new generic approach to combinatorial and multi-extremal optimization, as well as rare event simulation. It was motivated by an adaptative algorithm for estimating probabilities of rare events in complex stochastic networks, which involves variance minimization. In fact, it was soon realized that a simple cross-entropy modification could be used not only for estimating probabilities of rare events but for solving difficult combinatorial optimization problems as well. This is done by translating the deterministic optimization problem into a related stochastic optimization problem and then using rare event simulation techniques.\\
The naive idea to estimate rare events is to simulate huge samples of data. Another, less fastidious idea, should be to used Importance sampling, whose aim is to simulate the system according to a density, which should increase the occurence of this rare event. Whereas the determination of the tilting parameters used in the IS technique is quite hard, the CE method provides a way to cope efficently with such a phenomenom. \\
Let $S:\mathcal{X}\rightarrow \mathbb{R}$, be a real value function. 
We introduce $X=(X_1,X_2,...,X_N)$, which is a random vector defined on the space $\mathcal{X}$. 
Let $\{f(.;v)\}_v$ be a family of parametric densities with respect to the Lebesgue measure.\\
Actually, we want to estimate:

$$l\;=\;\mathbf{P}_{u}[\{S(X)\geq\gamma\}]\;=\;\mathbf{E}_u[\{S(X)\geq\gamma\}]\;.$$
 
If $l<10^{-5}$, we say that the event $\{S(X)\geq\gamma\}$, is a \textit{rare event}. Using IS, we try to simulate a random sample according to an importance sampling density g, on $\mathcal{X}$. As a result, we get an estimator of the form:
\begin{equation}
\hat{l}\;=\;\frac{1}{N}\displaystyle{\sum_{i=1}^{N}}\mathbf{1}_{\{S(X_i)\geq\gamma\}}\frac{f(X_i;u)}{g(X_i)}\;.
\end{equation}

The optimal zero variance associated estimator can easily be computed.
\begin{equation}
\label{IS}
g^{\star}(x)\;=\;\frac{\mathbf{1}_{\{S(x)\geq\gamma\}}f(x;u)}{l}\;.
\end{equation}

The idea in fact, is to choose g in the family of parametric densities $\{f(.;v)\}_v$, which is equivalent to determine the optimal associated parameter.
To determine this parameter, we will find the parametric density $f(.;v)$ which is the nearest from $g^{\star}$, using the Kullback-Leibler distance. This pseudo-distance between two densities $g$ and $h$, is defined as follows:
$$\mathcal{D}(g,h)\;=\;\mathbf{E}_{g}[\textrm{ln}\frac{g(X)}{h(X)}]\;=\;\int g(x)\;\textrm{ln}\;g(x)\;dx\;-\;\int g(x)\;\textrm{ln}\;h(x)\;dx\;.$$

As a result, minimizing the distance between $g^{\star}$ and $f(.;v)$ is equivalent to solving:
$$\max_{v}\int g^{\star}(x)\;\textrm{ln}\;f(x;v)\;dx\;.$$

By substitution of $(\ref{IS})$ into this equation, we get:

$$\max_v\;\mathcal{D}(v)\;=\;\max_v\mathbf{E}_u[\mathbf{1}_{\{S(X)\geq\gamma\}}\;\textrm{ln}\;f(X;v)]\;.$$

Finally, we can estimate v, using the associated stochastic problem:
\begin{equation}
v^{\star}\;=\;\arg\max_{v}\frac{1}{N}\sum_{i=1}^N\mathbf{1}_{\{S(X_i)\geq\gamma\}}\;\textrm{ln}\;f(X_i;v)\;.
\end{equation}

Now, we will try to highlight the link between rare event estimation and classical optimization problems. Consider an optimization problem of the form:
\begin{equation}
\label{opt_cont}
S(x^{\star})\;=\;\gamma^{\star}\;=\;\max_{x\in\mathcal{X}}S(x)\;.
\end{equation}

Our goal is to change this optimization problem into an estimation problem. 
Let $\{\mathbf{1}_{\{S(x)\geq\gamma\}},\;\gamma\in\mathbb{R}\}$, be a collection of indicator functions, and $\{f(.;v),v\in\mathcal{V}\}$, be a parametric family of densities.  
 
 For a fixed level, $u\in\mathcal{V}$, we associate to $(\ref{opt_cont})$, the following estimation problem:
 
 \begin{equation}
 l(\gamma)\;=\;\mathbf{P}_{u}[\{S(X)\geq\gamma\}]\;=\;\sum_{x}\mathbf{1}_{\{S(x)\geq\gamma\}}f(x;u)\;=\;\mathbf{E}_u[\mathbf{1}_{\{S(x)\geq\gamma\}}]\;.
 \end{equation}
 
 If $\gamma$ is close to $\gamma^{\star}$, then $f(.;v^{\star})$ will put the major part of its weight in $x^{\star}$. Consequently, the estimator developped in  the context of rare event simulation, can be put in aplication. However, to get a good estimator of that kind, it is necessary that $S(x)\geq\gamma$ for many realizations of the sample. This means that if $\gamma$ is close to $\gamma^{\star}$, then $u$ must be chosen so that $\mathbf{P}_u[\{S(x)\geq\gamma\}]$ remains not to small. The idea is to simultaneously simulate a sequence of levels $\hat{\gamma_1},\;\hat{\gamma_2},...,\hat{\gamma_T}$, and a sequence of parameters $\hat{v_1},\;\hat{v_2},...,\,\hat{v_T}$, such that $\hat{\gamma_T}$ tends towards the optimum $\gamma^{\star}$, and that $\hat{v_T}$ allows the density to give a higher weight to the states improving the performance. This bi-level algorithm takes the simple form:
 \begin{algo}
 $1$- Choose $\hat{v_0}\;=\;u$, and let $t=1$.\\
 $2$- Generate $X_1,\;...,\;X_N\;\sim\;f(.;v_{t-1})$, then compute the estimate of the $(1-\rho)$-quantile $\hat{\gamma_t}$ of the performance function.
 $$\hat{\gamma_t}\;=\;S_{([(1-\rho)N])}\;,$$
 where, $S_{([n])}\;,$is the n$^{th}$ element of the ordered statistics.\\
 $3$- Using $X_1,\;...,\;X_N$, solve the following stochastic problem:
 \begin{equation}
 \hat{v_t}\;=\;\arg\max_{v}\frac{1}{N}\sum_{i=1}^{N}\mathbf{1}_{\{S(X_i)\geq\hat{\gamma_t}\}}\;\textrm{ln}(f(X_i;v))\;.
 \end{equation}
 $4$- If $\hat{\gamma_t}=\hat{\gamma_{t-1}}=...=\hat{\gamma_{t-d}}$, STOP,\\
 else, set $t=t=1$.
 \end{algo}
 
 d is a constant ($d=5$ is generally a good compromise), and $\rho$ characterizes the level of rarity chosen.
 
 \subsubsection{Application of the CE method to estimate the amount of traffic on each VPN}
 
 We define a parametric utility function for each VPN. This utility function represents the subjective interpretation of the network manager on the impact of its bandwidth reservation on the quality of service for each VPN composing the whole MPLS network. 
 We suppose that the manager's interpretation follows a gamma density, whose parameter $p_i\in]1;\infty[,\;i=1,2,3$ is unknown. Indeed, it might be possible that the traffic sent on the VPN $1$, is not of the same type as the one on the VPN $2$. Hence, the operator would not use the same utilty function to characterize the impact of its allocation, on the traffic of the VPN $1$, or $2$. 
 
\begin{center}
\begin{figure}[h]
\includegraphics[scale=0.5]{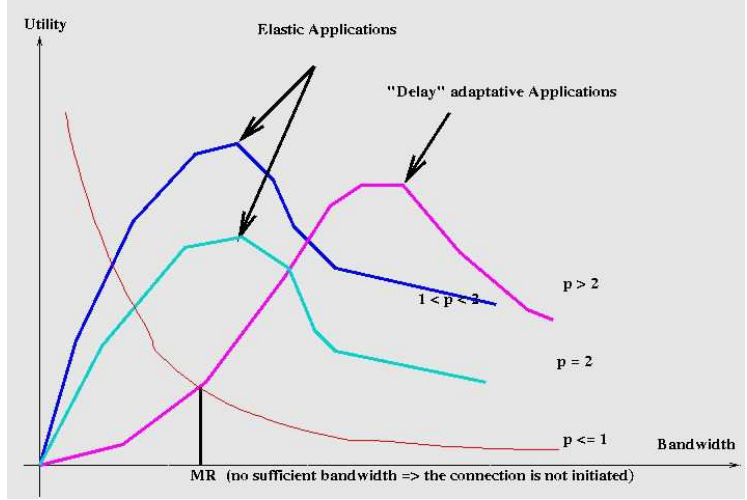}
\caption{The utility functions.}
\end{figure}
\end{center}

 The gamma density is particularly well adapted to model the two types of traffic, we have to deal with.
 The first one, associated with densities whose parameter $p\in]1;2]$, models \textit{elastic traffic}. This type of traffic have no real requirements in terms of delay and 
 transfer rates. Classical examples are email and data file transfer. The second one, for $p>2$, represents applications sensitive to delays, which requires an instantaneous transmission, like voice, or video over IP. If there is no sufficient bandwidth, the connection is not initialted, due to the existence of compression software which have an upper bound on the possible compression. This explains why the utility function equals zero below a certain value, MR (see $[17]$). 
 Furthermore, the fact that the densities decrease asymptotically, can be interpreted using an economical point of view. Indeed, if we reserve very large amounts of bandwidths for a particular traffic, the link capacity will be saturated, and the operator will make the client pay a lot, since there is no more ressource available for other type of traffic.
 
 We introduce a vector notation to store the amount of reserved bandwidth on each VPN link. For the VPN$1$ at the decision epoch $(t+1)$, we note:
 \begin{equation}
 B_X^{(t+1)}\;=\;\left( \begin{array}{c}
 B_{X_{12}}^{(t+1)}\\
 B_{X_{13}}^{(t+1)}\\
 B_{X_{21}}^{(t+1)}\\
 B_{X_{23}}^{(t+1)}\\ 
 B_{X_{31}}^{(t+1)}\\
 B_{X_{32}}^{(t+1)}
 \end{array} \right)\;.
 \end{equation}
 
 We proceed the same way to define $B_Y^{(t+1)}$ and $B_Z^{(t+1)}$, which are respectively the volumes of reserved bandwidth on the VPNs $2$ and $3$. 
 Besides, we suppose that these variables are generated following a gamma denstity. 
 
 \begin{equation}
 \left\{ \begin{array}{l}
 B_X^{(t+1)}\;\sim\;\gamma(p_1),\;p_1\in]1;\infty[,\\
 B_Y^{(t+1)}\;\sim\;\gamma(p_2),\;p_2\in]1;\infty[,\\
 B_Z^{(t+1)}\;\sim\;\gamma(p_3),\;p_3\in]1;\infty[.\\
 \end{array} \right.
 \end{equation}
 
 Actually, our aim, is to estimate the value of the unknown parameters $p_1,\;p_2$ and $p_3$, and to get a random sample solution of the equation $(\ref{MT})$. Recall that the gamma density $\gamma(p)\;p>0$, is of the form:
 $$f(x;p)\;=\;\frac{1}{\Gamma(p)}e^{-x}x^{p-1}\mathbf{1}_{\mathbb{R}^{+}}(x)\;,$$
 where, $\Gamma(p)\;=\;\int_{0}^{\infty}e^{-x}\;x^{p-1}\;dx\;.$
 
 \bigskip
 
 In terms of bandwidth reservation, we get the following matrix formulation:
 
 \begin{equation}
\label{MT}
(B^{\textrm{link}})^{(t+1)}\;=\;\left( \begin{array}{c}
B_1^{(t+1)}\\
B_2^{(t+1)}\\
\vdots\\
B_{|\mathcal{L}|}^{(t+1)}
\end{array} \right)\;=\;R\;\left( \begin{array}{c}
B_X^{(t+1)}\\
\hline\\
B_Y^{(t+1)}\\
\hline\\
B_Z^{(t+1)}
\end{array} \right)\;.
\end{equation}
 
On each link $(i,j)$ of the VPN$1$, knowing the volume of traffic $X_{ij}^{(t)}$ flowing through this link at the discrete time t, we can compute the minimum reserved bandwidth needed, solving a continuous optimization problem. 
The solution can be written using a bijective representation.
\begin{equation}
B_{X_{ij}}^{(t)}\;=\;\Phi(X_{ij}^{(t)})\;=\;X_{ij}^{(t)}+\frac{\sqrt{2\;X_{ij}^{(t)}}}{2\;p^X_{ij}}\;.
\end{equation}

\bigskip

\textbf{Hypothesis. We suppose that the sum of the parameters is constant:} $\mathbf{p_1+p_2+p_3\;=\;K,\;K>0\;.}$\\
This point takes into account the \textit{a priori} of the network manager, on the nature of the traffic he sent. For example, he may know that the traffic on the VPN$1$ and $2$ is elastic, which implies that $p_1+p_2\leq 4$, and fix an upper level on the density parameter for the VPN$3$, $2<p_3\leq 4$. Consequently, he gets the upper bound: $K=8$. 

\bigskip

Furthermore, since the random vectors $B_X^{(t+1)},\;B_Y^{(t+1)}$ and $B_Z^{(t+1)}$ are independent, the vector density is of the form: $\gamma(p)\;=\;\gamma(p_1)\;\gamma(p_2)\;\gamma(p_3)\;.$ In the application of the CE algorithm, we suppose that the decision epoch is fixed to $(t+1)$, but this algorithm can be applied anytime we have to jump from a global level to the local levels.

\bigskip

We will check the steps described in the CE algorithm.\\
$1$- We begin with the initialization of the density parameters.
$$\hat{p}^{(0)}\;=\;\left( \begin{array}{c}
\hat{p}_1^{(0)}\\
\hat{p}_2^{(0)}\\
\hat{p}_3^{(0)}
\end{array} \right)\;.$$
$2$- At the simulation instant $\tau\geq 1$, we simulate a sample of N vectors $B^1,\;B^2,...,\;B^N$, where $$B^i\;=\;\left( \begin{array}{c}
(B_X)^i\\
\hline\\
(B_Y)^i\\
\hline\\
(B_Z)^i
\end{array} \right)\;,\;i=1,2,...,N\;.$$

Hence, we infer the volumes of reserved bandwidths on each link.
\begin{equation}
(\hat{B^{\textrm{link}}})^i\;=\;
\left( \begin{array}{c}
B_1^{i}\\
B_2^{i}\\
\vdots\\
B_{|\mathcal{L}|}^{i}
\end{array} \right)\;=\;R\;\left( \begin{array}{c}
(B_X)^{i}\\
\hline\\
(B_Y)^{i}\\
\hline\\
(B_Z)^{i}
\end{array} \right)\;,\;i=1,2,...,N\;.
\end{equation}

Then we compute the performance function, 

$$S(B^i)\;=\;\frac{1}{\|(B^{\textrm{link}})^{(t+1)}-(\hat{B^{\textrm{link}}})^i\|}\;,\;i=1,2,..,N\;.$$

After ordering the statistic, we obtain the $(1-\rho)$-quantile of the performance function using the estimator:
$$\hat{\gamma_{\tau}}\;=\;S_{[(1-\rho)N]}\;.$$

$3$- Finally, to get the parameters updated, we have to solve a system of two equations.
\begin{equation}
\left\{ \begin{array}{l}
\frac{\textrm{ln}(\Gamma(p_1)\Gamma(K-p_1-p_2))}{p_1}\;=\;\frac{\displaystyle{\sum_{i=1}^N}\mathbf{1}_{\{S(B^i)\geq\hat{\gamma_{\tau}}\}}(\textrm{ln}((B^X_j)^i)-\textrm{ln}((B^Z_j)^i))}{\displaystyle{\sum_{i=1}^N}\mathbf{1}_{\{S(B^i)\geq\hat{\gamma_{\tau}}\}}}\;,\;\forall j=1,2,3\;,\\
\frac{\textrm{ln}(\Gamma(p_2)\Gamma(K-p_1-p_2))}{p_2}\;=\;\frac{\displaystyle{\sum_{i=1}^N}\mathbf{1}_{\{S(B^i)\geq\hat{\gamma_{\tau}}\}}(\textrm{ln}((B^Y_j)^i)-\textrm{ln}((B^Z_j)^i))}{\displaystyle{\sum_{i=1}^N}\mathbf{1}_{\{S(B^i)\geq\hat{\gamma_{\tau}}\}}}\;,\;\forall j=1,2,3\;.
\end{array} \right.
\end{equation}

If we restrict ourselves to integer values of the parameters $p_1$ and $p_2$, we just have to build a fine grid on the space defined by the equation $\{(p_1,p_2,p_3)|p_1+p_2+p_3=K,\;p_i\geq 0,\;i=1,2,3\}\;.$ We then get estimated values of the parameters. And, finally, we update the parameters to their new values:

$$p_1^{(\tau)}\;=\;p_1,\;p_2^{(\tau)}\;=\;p_2,\;p_3^{(\tau)}\;=\;p_3\;.$$

$4$- We stop as soon as: $\gamma_{\tau}\;=\;\gamma_{(\tau-1)}\;=\;...\;=\;\gamma_{(\tau-5)}\;.$

\bigskip

\textbf{Remark.} Once we have determined the optimal reserved bandwidth on each link, it is quite simple to get the value of traffic on the link, using the bijectivity of the function $\Phi$. 

%%%%%%%%%%%%%%%%%%%%%%%%%%%%%%%%%%%%%%%%%%%%%%%%%%%%%%%%%%%%%%
\begin{center}
\begin{figure}[h]
%\begin{tabular}{c}
\includegraphics[scale=0.6]{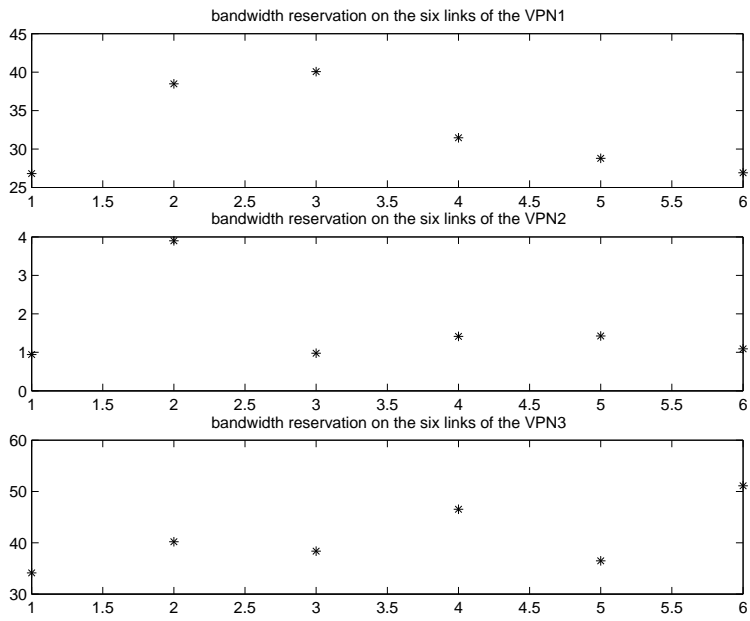}
\caption{Simulation of the reserved bandwidth, on each of the $6$ links of the $3$ VPNs, using the CE method. The constant was set to $K=70$, estimated values of the gamma densities parameters are: $p_1=3,\;p_2=4$ and $p_3=23$.}
\end{figure}
\end{center}
%%%%%%%%%%%%%%%%%%%%%%%%%%%%%%%%%%%%%%%%%%%%%%%%%%%%%%%%%%%%%%%
\begin{center}
\begin{figure}h]
\includegraphics[scale=0.5]{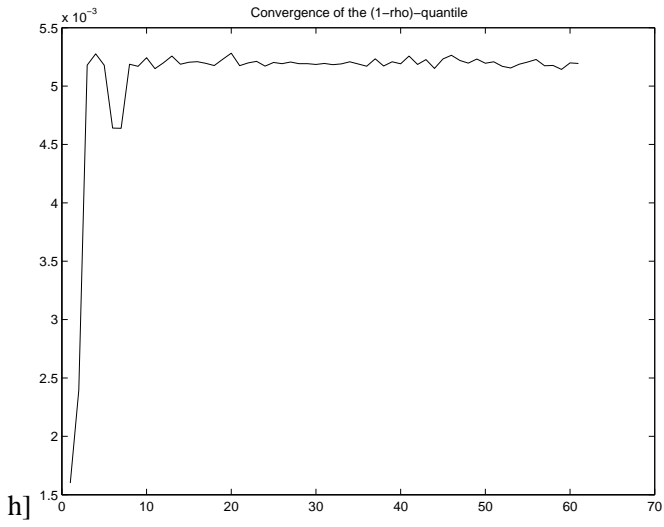}
\caption{The quantile function converges in around $15$ iterations.}
\end{figure}
%\end{tabular}
\end{center}
%%%%%%%%%%%%%%%%%%%%%%%%%%%%%%%%%%%%%%%%%%%%%%%%%%%%%%%%%%%%%%%

\subsection{Existence of stationary strategies for hierarchical MDPs}

The method we have developed so far, enables us to control \textit{optimaly} the dynamic evolution of our system, under the assumption that the horizon is finite. The optimality results from the introduction of centralized decisions, which aim to correct the evolution of the system, in order to satisfy the levels chosen by each VPN client. 

\begin{center}
\begin{figure}[h]
%\begin{tabular}{c}
\includegraphics[scale=0.7]{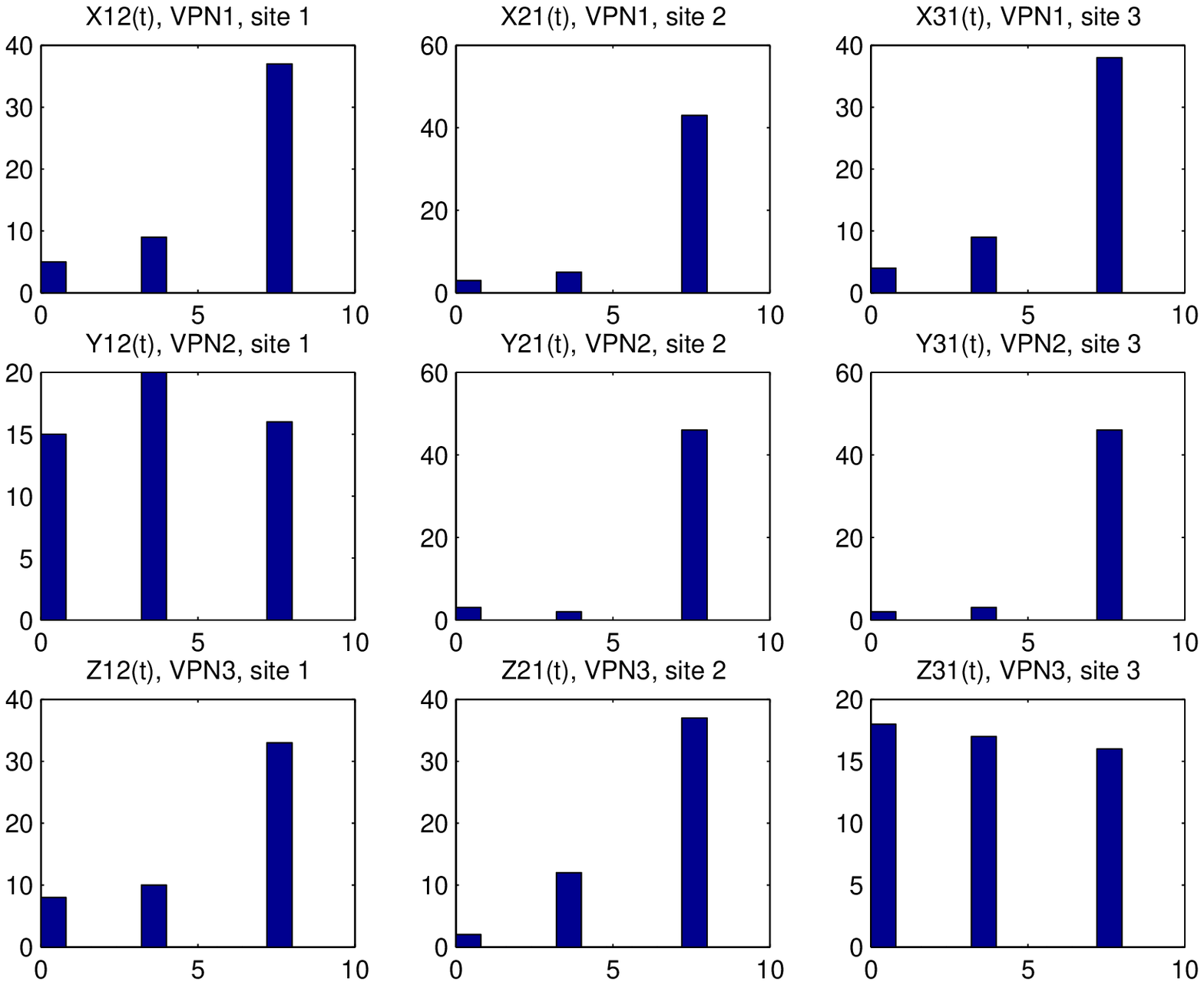}
\caption{\textbf{Hierarchical MDPs:} histograms of the dynamic evolution of the MDPS $X_1^{(t)},\;X_2^{(t)}$ and $X_3^{(t)}$. The state space associated to each MDP $X_i^{(t)},\;Y_i^{(t)}$ and $Z_i^{(t)},\;i=1,2,3$, is of cardinality $3$: $\mathbf{S}^{X_{i}}\;=\;\mathbf{S}^{Y{i}}\;=\;\mathbf{S}^{Z_{i}}\;=\;\{(0;9),(4;5),(8;1)\}\;,\;i=1,2,3$.
But, the global state space, which is required to take decisions at the global level, is of cardinality $27^3$. Consequently, it becomes fastly very hard to cope with such high dimensional spaces.}
\end{figure}
\end{center}

\begin{center}
\begin{figure}[h]
\includegraphics[scale=0.7]{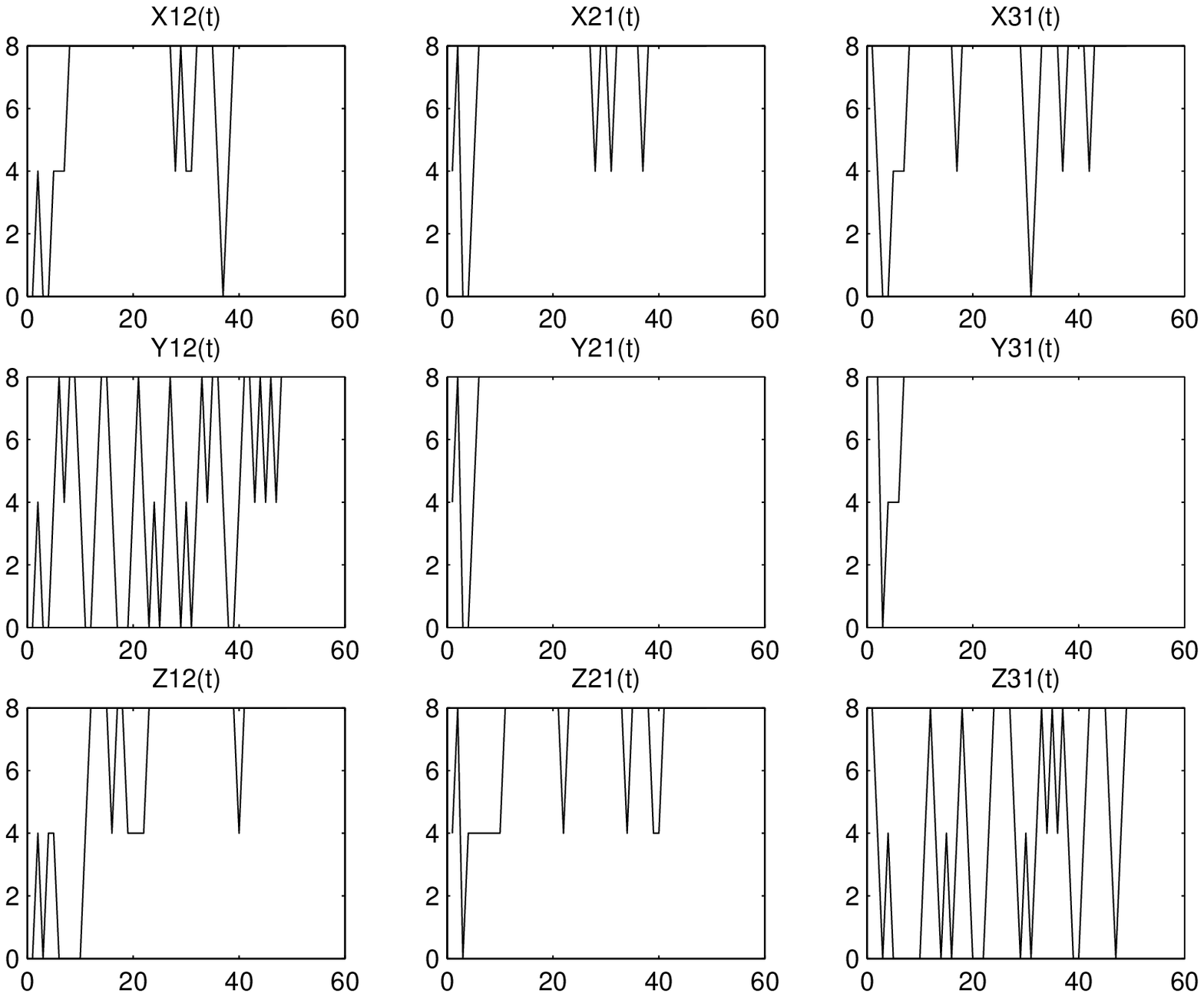}
\caption{\textbf{Hierarchical MDPs:} dynamic evolution of the traffic on each site, in each VPN. As an application, we consider a four state, state space. The hierarchical MDP approach enables us to control optimaly our system,
and furthermore, to characterize the evolution of each VPN comparatively to one another.}
\end{figure}
%\end{tabular}
\end{center}

In this section, our purpose is to study the asymptotic behavior of our system. The idea is, like in the very simple case of a $3$ site-VPN, to prove the existence of a stationary strategy via dynamic programming, and the convergence of the strategies obtained with the help of dynamic programming towards those stationary strategies. 
Indeed, we have proved the existence and the convergence of the strategies of a mono-path, stable routing VPN, towards a stationary strategy. 
Besides, the stationary strategies are deterministic, because we apply the simplex algorithm to compute them. Using the hierarchical MDP principle, Bellman's optimality equation gives us three distinct sequences of strategies, for each VPN, on the time interval $[0;T]$:

\begin{equation}
\left\{ \begin{array}{l}
(F_0^X,F_1^X,F_2^X,...,F_T^X)\rightarrow\;F^X\;,\;\textrm{for the VPN $1$},\\
(F_0^,F_1^Y,F_2^Y,...,F_T^Y)\rightarrow\;F^Y\;,\;\textrm{for the VPN $2$},\\
(F_0^Z,F_1^Z,F_2^Z,...,F_T^Z)\rightarrow\;F^Z\;,\;\textrm{for the VPN $3$}.
\end{array} \right.
\end{equation}

Where $F^X,\;F^Y$ and $F^Z$, are the associated asymptotic stationary strategies. 
If, we manage to prove that $(F_0^l,F_1^l,F_2^l,...,F_T^l)$ converges towards a stationary strategy $F^l$, for each link $l\in\mathcal{L}\;,$ 
then the system is asymptotically driven in a stable behavior, insofar as the local and the global strategies are both stationary. \\
Our aim presently, will be to prove that the stochastic process $\{L^{(t)}\}_t$, modeling the dynamic evolution of the demand on the link of the MPLS network, is ergodic. 
Consequently, for each link $l\in\mathcal{L}$, we need to solve a linear program of the form:

\begin{equation}
\left\{ \begin{array}{l}
\max\displaystyle{\sum_{s=1}^{\mathbf{S}^l}}\;\displaystyle{\sum_{a=1}^{3}}\mathcal{C}^{l}(s,a)\;x_{sa}^l\\
\displaystyle{\sum_{s=1}^{\mathbf{S}^l}}\;\displaystyle{\sum_{a=1}^{3}}[\delta(s,s')-\beta\;p(s'|s,a)]\;x_{sa}^l\;=\;\gamma(s),\;\forall s'\in\mathbf{S}^l\;,\\
x_{sa}^l\geq 0,\;a\in\mathbf{A},\;s\in\mathbf{S}\;.
\end{array} \right.
\end{equation}

For each state $s\in\mathbf{S}^l,\;l\in\mathcal{L}$, and each action $a\in\mathbf{A}_{\mathcal{L}}$, the optimal strategy on the link $l$, is easily obtained from:

\begin{equation}
f^l(s,a)\;=\;\frac{x_{sa}^l}{\displaystyle{\sum_{a\in\mathbf{A}}}x_{sa}^l}\;.
\end{equation}

All we need to do, is to verify the equality between temporal and spatial means:
\begin{equation}
\forall\;l\in\mathbf{S}^l,\;\forall k\in\mathbb{N},\;\lim_{t\rightarrow\infty}\sum_{s\in\mathbf{S}^{l}}s^k\;f_t^l(s,a)\;=\;\lim_{T\rightarrow\infty}\sum_{s\in\mathbf{S}^l}s^k\frac{\hat{x_{sa}^l}}{T}\;.
\end{equation}

Recall that the probability distributions $f^l(s,a)$ are obtained through linear programming, while we can infer the values of the parameters $\hat{x_{sa}^l}$ through simulation only.

\begin{center}
\begin{figure}[h]
\includegraphics[scale=0.5]{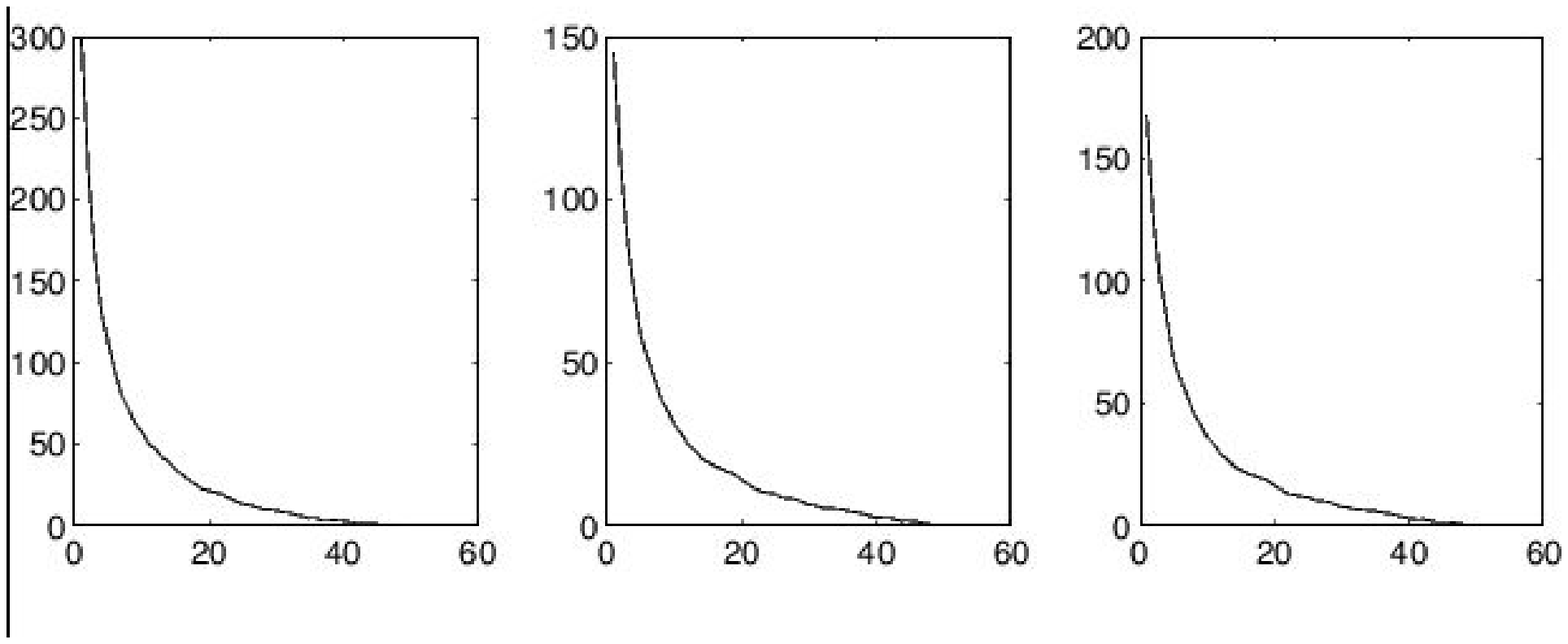}
\caption{\textbf{Hierarchical MDPs:} convergence of the empirical means of the MDP associated with the link $2$, $L_{2}^{(t)}$, conditionnaly to the choice of actions, $t\leq 300$.}
\end{figure}
\end{center}

We obtain the coincidence of these two means, which proves the ergodicity of the stochastic process $\{L^{(t)}\}_t$. 
We can conclude from these results that asymptotically, the local and the centralized strategies obtained via Bellman's optimality equation, will converge to stationary controls: $F^X,\;F^Y,\;F^Z$ and $F^L$, respectively. \\
We can infer from these results, that for $T$ not too large (i.e. $T\leq 200$), the stochastic dynamic approach is well adapted, but if choose to let $T$ increase towards infinity, it becomes rather tedious to compute all the optimal strategies. Since we have proved that asymptotically, our system adopts a stationary deterministic control, the use of linear programming provides an elegant and simple solution.

\subsection{The switching control game approach}

In this section, we delve into the fascinating world of stochastic Games. Our aim is still to determine stationary strategies. However, while we consider a hierarchical approach in the previous section, we try here, to model the problem as a matrix game, where decisions are alternatively taken either by the operator, either by the VPN owners, depending whether or not the satisfaction bounds are overwhelmed. Besides, the model is based on an initial assumption, which states that the game would converge towards an equilibrium, where the global delay on the links, and the sum of all the delays on each VPN would coincide, omitting an additive constant. Furthermore, switching-control games belong to the rare classes of games, which can be solved with the help of linear programming. Consequently, this model seems particularly promising.\\
We still consider an MPLS network, composed of $3$ independent VPNs. The routing is once more mono-path, and stable.\\
Using the vocabulary of Game Theory, we observe that our virtual network is made of various actors, whose interests are quite opposite. The only aim of the client, owning a VPN, is to get the best possible QoS for his personal traffic. The VPN owners behave completly non-cooperatively, since the traffic of each VPN  evolves without any collusion between the clients, whose single minded purposes are to minimize the delay on their own VPN. Furthermore, the clients are not aware of the presence of one another in the game, and behave perfectly selfishly.\\
As in the hierarchical case, we still assume that each VPN owner has previously determine a satisfaction level for its delay. In the case that this bound would be overwhelmed, the owners should have the opportunity to call for a centralized management. The operator realizes this centralized management, by choosing controls on the links of the network. The global traffic is still supposed to follow the worst possible evolution, but the operator's purpose is now to minimize the global delay. \\

Once more, we refer to the global process, $(X^{(t)},Y^{(t)},Z^{(t)})$, which takes its values in the global state space $\mathbf{S}^X\times\mathbf{S}^Y\times\mathbf{S}^Z\;.$

If we think about the way our decisions are made, we realize that the state space can be partionned into two disjoint subsets. 
Indeed, there are some combinations of states $(x^{(t)},y^{(t)},z^{(t)})$, that will automatically violate the satisfaction bounds imposed by at least one of the owner. As a result, for all these global states, the network should be centrally controled. This subset of the state space will be denoted $\mathbf{E}^2$.\\
On the contrary, on the rest of the global states called $\mathbf{E}^1$, the satisfaction bounds are not overwhelmed, and the decisions are taken independently on each VPN.

\begin{center}
\begin{figure}[h]
\includegraphics[scale=0.4]{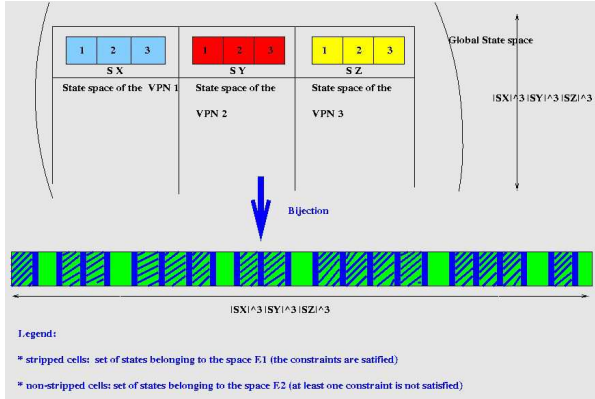}
\caption{The classification of the state space of the $3$ VPNs.}
\end{figure}
\end{center}

The idea is to model the problem as a two-person zero-sum game. The first player will represent the set of the $3$ VPNs, evolving independently and selfishly. The action space $\mathbf{A}_{\mathcal{G}}$ contains all the possible combinations for the choices of each site, in each VPN. Naturally each possible combination is formally represented as a $9$-dimensional vector. 
We assume that the choice of action on each VPN link is reduced to the $3$ alternatives described in the section $2$: $\{a_0;a_1;a_2\}\;.$ \\
Consequently, $\mathbf{A}_{\mathcal{G}}$ is of finite cardinality, since the choice of actions on the $3$ VPNs is independent, and that the action space for each VPN is finite.\\
The second player will stand for the network manager, who should centrally manage the whole MPLS network, by taking actions on the links of the network. This time, the action space is denoted $\mathbf{D}$, and it contains all the possible combinations of actions that could be chosen on each link. We suppose that each decision on each link, is chosen in the $3$ element space: $\{d_0;d_1;d_2\}$. It is allowed to have the same choices of basic actions: $d_i\;=\;a_i,\;\forall i$.\\
We make the assumption which characterizes a switching-control game: \textit{on the states belonging to $\mathbf{E}^1$, only player $1$ can influence the transitions, whereas on the states belonging to $\mathbf{E}^2$, it is the player $2$ who controls the transitions.}\\
However, the reward function depends on the actions of both players, and takes the formal form:

\begin{eqnarray}
r(s,A,D)\;&=&\;\sum_{i,j\in\{1,2,3\},i\neq j}\mathcal{C}^X(s_{X_{ij}},a_{X_{ij}})+\sum_{i,j\in\{1,2,3\},i\neq j}\mathcal{C}^Y(s_{Y_{ij}},a_{Y_{ij}}){}\nonumber\\
& & {}+\sum_{i,j\in\{1,2,3\},i\neq j}\mathcal{C}^Z(s_{Z_{ij}},a_{Z_{ij}})-\left[\sum_{i\in\mathcal{L}}\mathcal{C}^L((R\;(s_X|s_Y|s_Z))_i,d_i)+\lambda\right]\;,{}\nonumber\\
& & {}\forall s\in\mathbf{E}^1\cup\mathbf{E}^2,\;\forall A\in\mathbf{A},\;\forall D\in\mathbf{D},\;\lambda\in\mathbb{R}^{+}\;.
\end{eqnarray}

$\lambda\in\mathbb{R}^{+}$, should model the amount of bandwith that the operator always keeps free in the fear of congestion.\\
$s\;=\;\left( \begin{array}{l l l l l l l l l l l l}
s_{X_{12}} s_{X_{21}} s_{X_{31}} | s_{Y_{12}} s_{Y_{21}} s_{Y_{31}} | s_{Z_{12}} s_{Z_{21}} s_{Z_{31}}
\end{array} \right)^{T}\in\mathbf{E}^1\cup\mathbf{E}^2$, represents a realization of the global process taking value in the state space. For the ease of notations, we would rather use the following one:
$s\;=\;(s_X|s_Y|s_Z)^T\;.$\\
Furthermore, to each state belonging to $\mathbf{E}^{2}$, we can define a specific configuration for the traffic value on the MPLS links. $R$, being the routing matrix, we obtain the values of the traffic on the links, by computing the matrix equation:

$$l\;=\;R\; \left( \begin{array}{l}
s_{X_{12}}\\
t_{out1}^X-s_{X_{12}}\\
s_{X_{21}}\\
t_{out2}^X-s_{X_{21}}\\\
s_{X_{31}}\\
t_{out3}^X-s_{X_{31}}\\
\hline\\
s_{Y_{12}}\\
t_{out1}^Y-s_{Y_{12}}\\
s_{Y_{21}}\\
t_{out2}^Y-s_{Y_{21}}\\\
s_{Y_{31}}\\
t_{out3}^Y-s_{Y_{31}}\\
\hline\\
s_{Z_{12}}\\
t_{out1}^Z-s_{Z_{12}}\\
s_{Z_{21}}\\
t_{out2}^Z-s_{Z_{21}}\\\
s_{Z_{31}}\\
t_{out3}^Z-s_{Z_{31}}\\
\end{array} \right)\;:=\;R\;\left( \begin{array}{l}
\tilde{s_X}\\
\hline\\
\tilde{s_Y}\\
\hline\\
\tilde{s_Z}\\
\end{array} \right)\;.$$

The set of the various possible link traffic configurations, will be once more, denoted $\mathcal{L}$.\\

The resulting game is a zero-sum game, since player $1$ wants clearly to maximize the reward function, whereas player $2$ wants to minimize the reward.\\
We know, as is the case for every matrix game, that the game has a value, and that both players do possess optimal actions. This famous result is due to J. von Neumann, and can be found in the rich literature related to the subject. The value of the game at time t, will be denoted $V_t$. It is a vector of length $|\mathbf{E}^1|+|\mathbf{E}^2|$.\\ 
Presently, our purpose will be to determine the optimal stationary strategies associated to each player. To this end, we use the algorithm developped in $[5]$, which is proved to converge in a finite number of iterations.\\
In fact, switching-control games belong to the class of stochastic games satisfying the order field property, which characterizes the single class of games, whose solution can be found in the same algebraic field as the data of the game. This class of games is all the more important that, only for such games, can one expect to be able to develop finite algorithm for deriving a solution.\\
Now, we will describe in details the algorithm we have used, to compute the stationary strategies.

\bigskip

 $\star$ We start by choosing an initial deterministic strategy for the player $1$, that will be noted $F_0(\mathbf{E}^1)$. This formal presentation, only means that we choose a pure action in each state of the state space $\mathbf{E}^1$.
 
 \bigskip
 
 $\star$ Then, $F_t(\mathbf{E}^1)$ being fixed, we solve the discounted game with a single controller, $\beta\in[0;1[$:

\begin{equation}
\left\{ \begin{array}{l}
 \max\displaystyle{\sum_{s=1}^{|\mathbf{E}^2|}}\;\displaystyle{\sum_{d\in\mathbf{D}}}\;\{\sum_{i\in\mathcal{L}}[\mathcal{C}^L((R\;\left( \begin{array}{l}
\tilde{s_X}\\
\hline\\
\tilde{s_Y}\\
\hline\\
\tilde{s_Z}\\
\end{array} \right))_i,d)+\lambda]\;x_{sd}\}\\
 \displaystyle{\sum_{s=1}^{|\mathbf{E}^2|}}\;\displaystyle{\sum_{d\in\mathbf{D}}}[\delta(s,s')-\beta\;p(s'|s,d)]\;x_{sd}\;=\;\gamma(s),\;\forall\;s'\in\mathbf{S}\\
 x_{sd}\geq 0,\;\forall d\in\mathbf{D},\;\forall s\in\mathbf{E}^2\;.
 \end{array} \right.
 \end{equation}
 
 In the value vector $V_t$, we stock the value of the game for each component belonging to $\mathbf{E}^2$. 
 
 \bigskip
 
 $\star$ For each state $s\in\mathbf{E}^1$, we will determine the action $F_{t+1}(s)$ as an extreme optimal action for player $1$, in the matrix game:
 \begin{equation}
 \mathcal{R}_{\beta}(s,V_t)\;=\;\left[(1-\beta)\;r(s,a,d)+\beta\sum_{s'=1}^{[\mathbf{E}^1|+|\mathbf{E}^2|}p(s'|s,a)V_t(s')\right]\;,\;\forall a\in\mathbf{A},\;\forall d\in\mathbf{D}\;.
 \end{equation}
 
 \bigskip
 
$\star$ If $V_t=V_{t-1}$, then $V_t$ is the value of the game , and $F_{t}(\mathbf{E}^1)$ is the projection of an optimal stationary strategy for the game on the space $\mathbf{E}^1$.

\begin{center}
\begin{figure}[h]
\includegraphics[scale=0.7]{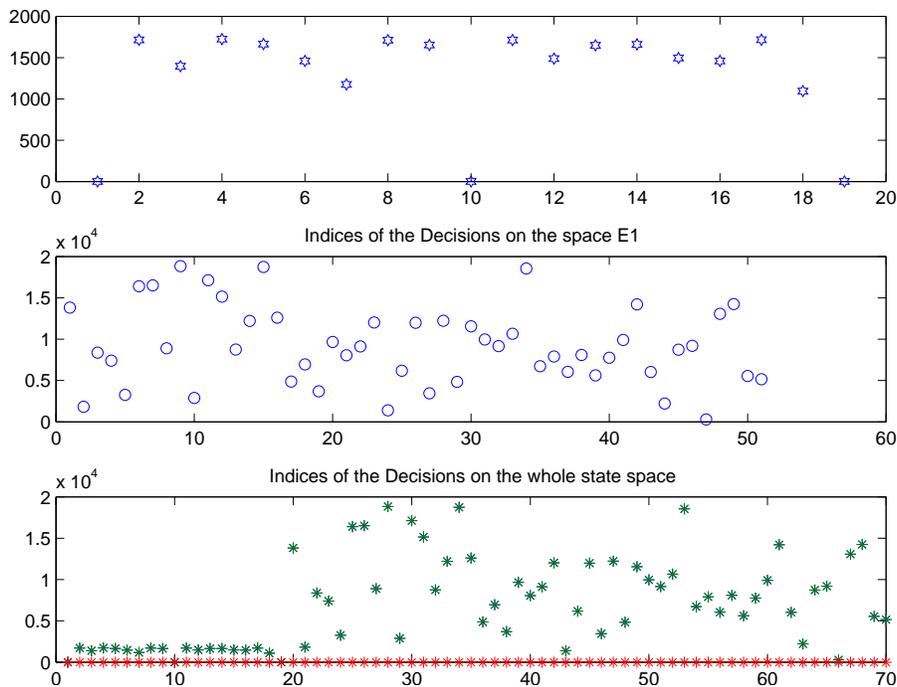}
\caption{The algorithm gives us the indices of the optimal action associated to each state. The cardinality of the action space associated to the $3$ VPNs was of $27^3$, and the dimension of the action space associated to the links was reduced here, to $12^3$. On the first picture, we can see the indices of the optimal actions on the space $\mathbf{E}^2$ of the global decisions, on the second one, we have represented the indices of the local strategies, on the space $\mathbf{E}^1$.}
\end{figure}
\end{center}

\bigskip

\textbf{Remark.} Once we have determined the stationary strategies, we use the same principle as in the hierarchical approach to construct the traffic trajectories. We begin to choose an initial state. Then, if this state belongs to the subspace $\mathbf{E}^1$, we take the associated optimal actions on each VPN, but if the state belongs to $\mathbf{E}^2$, we choose the optimal global actions on the links. One of the advantage of this approach, is that we do not need to verify whether or not, the satsfaction levels are overwhelmed, since the localisation of the state on the state space gives us enough information.

\section{The curse of dimensionality and optimization in policy space}
The use of Markov decision processes and the associated dynamic programming methodology become rapidly limited, due to the high cardinality of the state space. A solution to such a problem, lies in the introduction of parametric representations. There are three main methods to tackle the problem of dimensionality. The first well-known method called neuro-dynamic programming, or reinforcement learning, requires the introduction of weights in the value function. In each state $s\in\mathbf{S}$, the value function takes the form: $V(s,r),\;r\geq0$. The idea is to tune the weights, so as to obtain a good approximation of the value function, and to infer a policy as close as possible to the optimal one. The second method, essentialy developed in $[18]$, considers a class of policies described by a parameter vector $\theta\in\mathbb{R}^K$. The policy is improved by updating $\theta$ in a gradient direction, via simulation. The third and last one, called actor-critic, combines the principles of both approaches.\\
In this article, we concentrate our study on the improvement of the parametrized policy through the policy space. 
Our performance metric will be the average reward function, since the methodology developed in $[18]$ requires such an assumption in order to introduce the steady state probabilities in the performance function and later, to derive a proper estimate for the gradient function. The long term average reward is commonly denoted:
\begin{equation}
\lambda(\theta)\;=\;\lim_{T\rightarrow\infty}\frac{1}{T}\mathbf{E}\left[\sum_{t=0}^T\mathcal{C}_t(X^{(t)},\theta)|X^{(0)}=s\right]\;,
\end{equation} 
where, we still have, 
$$\mathcal{C}_t(X^{(t)},\theta)\;=\;\sum_{i,j}\left[\frac{X_{ij}^{(t)}}{\Phi(X_{ij}^{(t)})-X_{ij}^{(t)}}+p_{ij}\left(\Phi(X_{ij}^{(t)})-\Phi(X_{ij}^{(t-1)})\right)\right]\;.$$

At the instant t, we define a parametric matrix of strategies $F_t(\theta)$. Let $\theta\in\mathbf{R}^{K}$, be a parameter of size $K>0$. We define $f_t(s,a,\theta)$, as the probability to be in the state $s\in\mathbf{S}$, while we choose the action $a\in\mathbf{A}$, at the decision  epoch t.  
The parametrized transition probabilities and reward function, take the form:

\begin{equation}
\label{ptheta}
\left\{ \begin{array}{l}
p_{\theta}(s,s')\;=\;\displaystyle{\sum_{a\in\mathbf{A}}}f_t(s,a,\theta)\;p(s'|s,a)\;,\;\forall s,s'\in\mathbf{S}\;,\\
\mathcal{C}_t(s,\theta)\;=\;\displaystyle{\sum_{a\in\mathbf{A}}}f_t(s,a,\theta)\;\mathcal{C}_t(s,a)\;,\;\forall s\in\mathbf{S}\;.
\end{array} \right.
\end{equation}

We denote, $\mathcal{P}\;=\;\{P(\theta)\;=\;(p_{\theta}(s,s'))_{s,s'\in\mathbf{S}},\;\theta\in\mathbf{R}^{K}\}$, the set of the transition probabilities, and $\bar{\mathcal{P}}$, its closure which is also composed of stochastic matrices.   
Furthermore, we make the following assumption, required to prove the convergence of the associated algorithm:\\
$\bullet$ \textit{The Markov chain corresponding to every $P\in\bar{\mathcal{P}}$, is aperiodic, which means that the GCD of the length of all its cycles is one. Besides, there exists a state $s^{\star}\in\mathbf{S}$, which is recurrent for every such Markov chain.}\\

Our purpose is presently, to maximize the average reward:
\begin{equation}
\theta^{\star}\;=\;\displaystyle{\arg\max_{\theta}}\left\{\lambda(\theta)\right\}\;=\;\displaystyle{\arg\max_{\theta}}\left\{\lim_{T\rightarrow\infty}\frac{1}{T}\mathbf{E}\left[\sum_{t=0}^T\mathcal{C}_t(x^{(t)},\theta)|x^{(0)}=s\right]\right\}\;,
\end{equation}

where $x^{(t)}$, represents a realization of the stochastic process $\{X^{(t)}\}_t$, at the instant t. The esperance is computed relatively to the randomized strategy $F_t(\theta)$. 
The first idea is to introduce the well-known gradient algorithm, to get an estimate of the parameter. 
$$\theta(t+1)\;=\;\theta(t)\;+\;\gamma_t\;\nabla_{\theta}\;\lambda\left(\theta(t)\right),\;\gamma_t\;=\;\frac{1}{t},\;t\in\mathbb{N}\;.$$
Unfortunately, we can't compute analytically the gradient of the performance function, and must resort to simulation.  
The algorithm developed in $[18]$, updates at every time step the value of the parameter, and uses a biased estimate (whose bias asymptotically vanishes) of the gradient of the performance metric. 

\begin{equation}
\left\{ \begin{array}{l}
\theta(t+1)\;=\;\theta(t)+\gamma_t\;\left(\nabla_{\theta}\;\mathcal{C}_t(x^{(t)},\theta(t))\;+\;\left(\mathcal{C}_t(x^{(t)},\theta(t))-\tilde{\lambda}\right)\;z_t\right)\;,\\
\tilde{\lambda_{t+1}}\;=\;\tilde{\lambda_t}\;+\;\eta\;\gamma_t\;\left(\mathcal{C}_t(x^{(t)},\theta(t))-\tilde{\lambda_t}\right)\;.
\end{array}\right.
\end{equation}

$\eta>0$, is a parameter which enables us to scale the stepsize of our algorithm for updating $\tilde{\lambda_t}$ by a positive constant. 
Then, we simulate a transition to the next state $x^{(t+1)}$ following the transition probabilities $\{p_{\theta_{t+1}}(x^{(t)},s)\;,\;s\in\mathbf{S}\}\;.$

At the same time, $z$ is updated according to the following rules:

\begin{equation}
z_{t+1}\;=\;\left\{ \begin{array}{l}
0,\;\textrm{if}\;x^{(t+1)}=s^{\star},\\
z_t+L_{\theta_t}(x^{(t)},x^{(t+1)})\;,\;\textrm{otherwise}\;,
\end{array} \right.
\end{equation}

where $L_{\theta_t}(x^{(t)},x^{(t+1)})\;=\;\frac{\nabla_{\theta}p_{\theta_t}(x^{(t)},x^{(t+1)})}{p_{\theta_t}(x^{(t)},x^{(t+1)})}$, if $p_{\theta_t}(x^{(t)},x^{(t+1)})>0$, $0$ otherwise.
This term can be interpreted as a likelihood ratio derivative term.

In the case of a $3$ site-VPN, we choose simple parametric strategies. For the site $i$, $i\in\{1,2,3\}$, we set:

\begin{equation}
\left\{ \begin{array}{l}
(f_t^1)^{\textrm{site $i$}}\;=\;f_t^{\textrm{site $i$}}(X_i^{(t)},a_0,\theta)\;=\;\frac{1}{1+\exp\left[\left(\Phi(X_{ij}^{(t)})+\Phi(X_{ik}^{(t)})\right)-\theta_1^{\textrm{site $i$}}\right]},\\
\textrm{the probability to choose the action $a_0$ for the chain $\{X_{i}^{(t)}\}_t$},\;k\in\{1,2,3\},i\neq\;j,\;i\neq\;k\;,\\
(f_t^2)^{\textrm{site $i$}}\;=\;f_t^{\textrm{site $i$}}(X_i^{(t)},a_1,\theta)\;=\;\frac{1}{1+\exp\left[\left(\Phi(X_{ij}^{(t)})+\Phi(X_{ik}^{(t)})\right)-\theta_2^{\textrm{site $i$}}\right]},\\
\textrm{the probability to choose the action $a_1$ for the chain $\{X_{i}^{(t)}\}_t$},\\
(f_t^3)^{\textrm{site $i$}}\;=\;f_t^{\textrm{site $i$}}(X_i^{(t)},a_2,\theta)\;=\;\frac{1}{1+\exp\left[\left(\Phi(X_{ij}^{(t)})+\Phi(X_{ik}^{(t)})\right)-\theta_3^{\textrm{site $i$}}\right]},\\
\textrm{the probability to choose the action $a_2$ for the chain $\{X_{i}^{(t)}\}_t$}.
\end{array} \right.
\end{equation}

We note that:
 \begin{equation}
\left\{ \begin{array}{l}
f_t^1\geq 0.5 \Leftrightarrow \left(\Phi(X_{ij}^{(t)})+\Phi(X_{ik}^{(t)})\right)\leq \theta_1\;,\\
f_t^2\geq 0.5 \Leftrightarrow \left(\Phi(X_{ij}^{(t)})+\Phi(X_{ik}^{(t)})\right)\leq \theta_2\;,\\
f_t^3\geq 0.5 \Leftrightarrow \left(\Phi(X_{ij}^{(t)})+\Phi(X_{ik}^{(t)})\right)\leq \theta_3\;.\\
\end{array} \right.
\end{equation}
 
As a result, the parameters $\theta_i,\;i=1,2,3$ can be interpreted as fuzzy bounds for the system, since it determines the probability to choose the action i.

\begin{center}
\begin{figure}[h]
\includegraphics[scale=0.4]{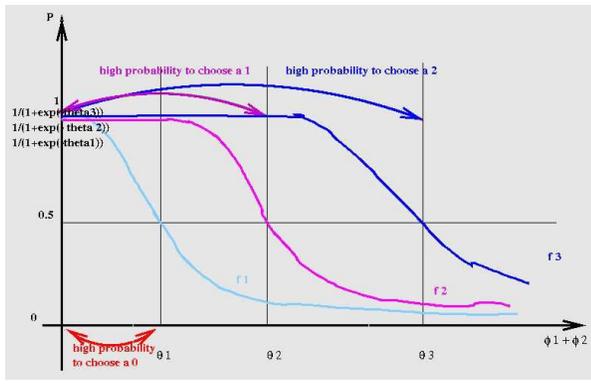}
\caption{Belief functions, or parametrized strategies.}
\end{figure}
\end{center}

Remind that the cost function for each site i ($i=1,2,3$), of the VPN$1$, is of the form:

\begin{equation}
\label{cost_par}
\mathcal{C}_t(X_i^{(t)})\;=\;\frac{X_{ij}^{(t)}}{\Phi(X_{ij}^{(t)})-X_{ij}^{(t)}}+p_{ij}\;(\Phi(X_{ij}^{(t)})-\Phi(X_{ij}^{(t-1)}))\;+\;\frac{X_{ik}^{(t)}}{\Phi(X_{ik}^{(t)})-X_{ik}^{(t)}}+p_{ik}\;(\Phi(X_{ik}^{(t)})-\Phi(X_{ik}^{(t-1)}))\;,
\end{equation}

where $j,k\in\{1,2,3\},i\neq j,\;i\neq k,\;j\neq k\;.$ \\
Besides, we infer the analytical expression of $\mathcal{C}_t(X_i^{(t)},\theta(t))$ from $(\ref{ptheta})$ and $(\ref{cost_par})$.\\
 
Finally, the iterative algorithm applied to the site $j$ $(j=1,2,3)$ of the VPN$1$, takes the form:

\begin{equation}
\label{tsitsi}
\left\{ \begin{array}{l}
\theta_1^{\textrm{site $j$}}(t+1)\;=\;\theta_1^{\textrm{site $j$}}(t)+\gamma_t\;\left(\nabla_{\theta_1^{\textrm{site $j$}}}\;\mathcal{C}_t(x^{(t)},\theta(t))+\left(\mathcal{C}_t(x^{(t)},\theta(t))-\tilde{\lambda_t^{\textrm{site $j$}}}\right)z_1^{\textrm{site $j$}}(t)\right)\;,\\
\theta_2^{\textrm{site $j$}}(t+1)\;=\;\theta_2^{\textrm{site $j$}}(t)+\gamma_t\;\left(\nabla_{\theta_2^{\textrm{site $j$}}}\;\mathcal{C}_t(x^{(t)},\theta(t))+\left(\mathcal{C}_t(x^{(t)},\theta(t))-\tilde{\lambda_t^{\textrm{site $j$}}}\right)z_2^{\textrm{site $j$}}(t)\right)\;,\\
\theta_3^{\textrm{site $j$}}(t+1)\;=\;\theta_3^{\textrm{site $j$}}(t)+\gamma_t\;\left(\nabla_{\theta_3^{\textrm{site $j$}}}\;\mathcal{C}_t(x^{(t)},\theta(t))+\left(\mathcal{C}_t(x^{(t)},\theta(t))-\tilde{\lambda_t^{\textrm{site $j$}}}\right)z_3^{\textrm{site $j$}}(t)\right)\;,\\
\tilde{\lambda_{t+1}^{\textrm{site $j$}}}\;=\;\tilde{\lambda_t^{\textrm{site $j$}}}+\eta\;\gamma_t\;\left(\mathcal{C}_t(x^{(t)},\theta(t))-\tilde{\lambda_t^{\textrm{site $j$}}}\right)\;,
\end{array} \right.
\end{equation}

where $x^{(t)}\;=\;(x_{1}^{(t)},x_{2}^{(t)})$, is a realization of the process $X_j^{(t)}\;=\;(X_{ji}^{(t)},X_{jk}^{(t)})$ $(j=1,2,3)$, and,

\begin{equation}
z_i^{\textrm{site $j$}}(t+1)\;=\;\left\{ \begin{array}{l}
0,\;\textrm{if}\;x^{(t+1)}\;=\;s^{\star}\;,\\
z_i^{\textrm{site $j$}}(t)\;+\;L_{\theta_i^{\textrm{site $j$}}(t)}(x^{(t)},x^{(t+1)}),\;\textrm{otherwise}\;(i=1,2,3)\;.
\end{array} \right.
\end{equation}

%%%%%%%%%%%%%%%%%%%%%%%%%%%%%%%%%%%%%%%%%%% 
\begin{figure}[h]
\includegraphics[scale=0.45]{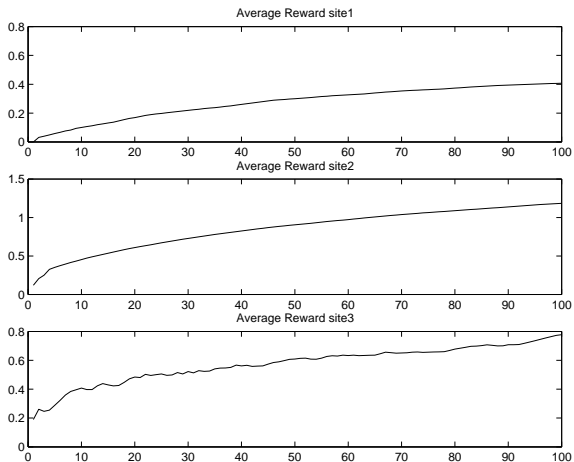}
\caption{\textbf{Parametric strategies on a single VPN.} As an example, we consider once more, the $3$-site VPN, with stable and mono-path routing
(i.e. the unique path between each couple of nodes, is the directed link between these sites). The sites are independent of one another, consequently the algorithm applies independently. The recurrent states are supposed to be: $(t_{\textrm{out}}^{1};0)$ for the site $1$, $(t_{\textrm{out}}^{2};0)$ for the site $2$, and $(t_{\textrm{out}}^{3};0)$, for the third site. The convergence of the average reward occures in around $100$ iterations.}
\end{figure}

\begin{figure}[h]
\includegraphics[scale=0.45]{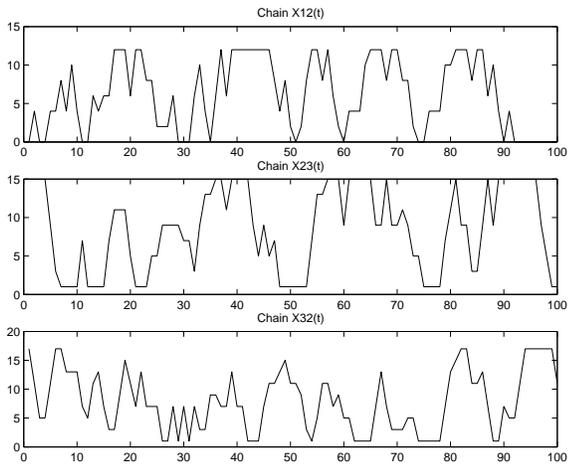}
\caption{\textbf{Parametric strategies on a single VPN.} The second picture represents the dynamic evolution of the sampled trajectories for each site of the VPN.}
\end{figure}

\begin{figure}[h]
\includegraphics[scale=0.45]{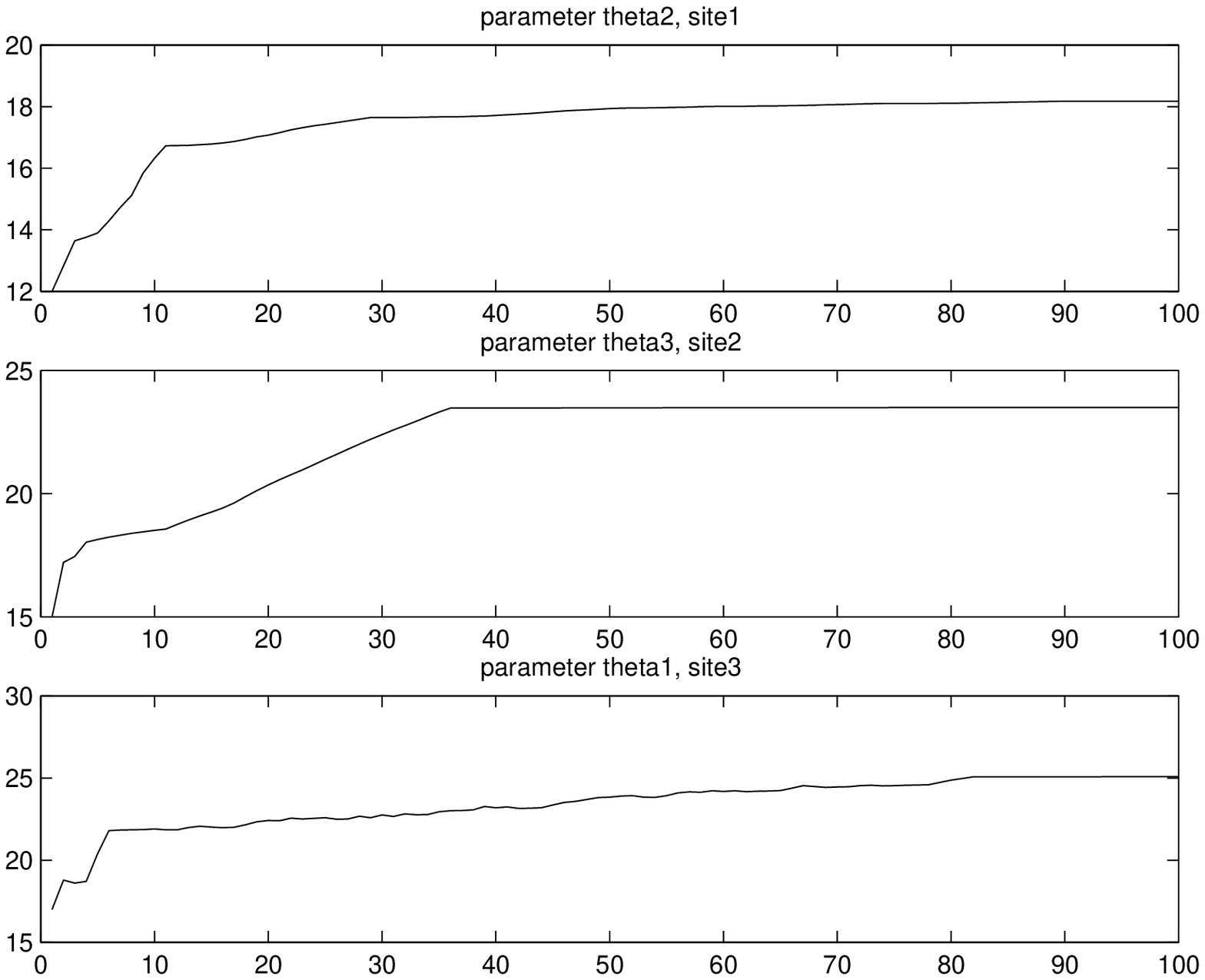}
\caption{\textbf{Parametric strategies on a single VPN.} Convergence of the parameters defining the parametric strategies of the VPN$1$'s sites.}
\end{figure}
%%%%%%%%%%%%%%%%%%%%%%%%%%%%%%%%%%%%%%%%%%%%

\bigskip

The application of the algorithm to a MPLS network of $3$ independent VPNs, with a stable and mono-path routing, and the possiblity to introduce a central management, requires the estimation of $3*12=36$ additional parameters. Locally we still have to cope with $3$ independent systems of the form $(\ref{tsitsi})$, and consequently we generate $3$ samples per VPN, at the decision epoch t: $(x_1^{(t)},x_2^{(t)},x_3^{(t)})\sim\;p_{\theta_t^{\textrm{site $1$}}}(s,s'),\;s,\;s'\in\mathbf{S}^1$ for the VPN$1$, $(y_1^{(t)},y_2^{(t)},y_3^{(t)})\sim\;p_{\theta_t^{\textrm{site $2$}}}(s,s'),\;s,\;s'\in\mathbf{S}^2$ for the VPN$2$, and finally, $(z_1^{(t)},z_2^{(t)},z_3^{(t)})\sim\;p_{\theta_t^{\textrm{site $3$}}}(s,s'),\;s,\;s'\in\mathbf{S}^3$ for the VPN$3$.

We introduce global transition probabilities, and rewards on the links:

\begin{equation}
\left\{ \begin{array}{l}
p_{\theta^{\textrm{link}}_t}(l,l')\;=\;\displaystyle{\sum_{d\in\mathbf{D}}}f_t(l,d,\theta^{\textrm{link}}_t)\;p(l'|l,d)\;,\\
\mathcal{C}_t(l,\theta)\;=\;\displaystyle{\sum_{d\in\mathbf{D}}}f_t(l,d,\theta)\;\mathcal{C}_t(l,d)\;,\;l,\;l'\in\mathbf{S}^{1}\cup\mathbf{S}^{2}\cup...\cup\mathbf{S}^{|\mathcal{L}|}\;.
\end{array} \right.
\end{equation}    

Each time, one of the satisfaction levels is overwhelmed, we solve the global iterative algorithm:

\begin{equation}
\left\{ \begin{array}{l}
\theta_{t+1}^{\textrm{link}}\;=\;\theta_{t}^{\textrm{link}}\;+\;\gamma_t\;\left[\nabla_{\theta}\;\mathcal{C}_t\left(R(x^{(t)},y^{(t)},z^{(t)})^T,\theta_t\right)\;+\;\left(\mathcal{C}_t\left(R(x^{(t)},y^{(t)},z^{(t)})^T,\theta_t\right)-\tilde{\lambda_t^{\textrm{link}}}\right)\;z_t\right]\;,\\
\tilde{\lambda_{t+1}^{\textrm{link}}}\;=\;\tilde{\lambda_t^{\textrm{link}}}+\eta\;\gamma_t\;\left[\mathcal{C}_t\left(R(x^{(t)},y^{(t)},z^{(t)})^T,\theta_t\right)\;-\;\tilde{\lambda_t^{\textrm{link}}}\right]\;,
\end{array} \right.
\end{equation}

with,

$$z_{t+1}\;=\;\left\{ \begin{array}{l}
0,\;\textrm{if}\;l^{(t+1)}\;=\;R\;[t_{\textrm{out}}^1\;0\;|\;t_{\textrm{out}}^2\;0\;|\;t_{\textrm{out}}^3\;0]^T\;,\\
z_t\;+\;L_{\theta_t^{\textrm{link}}}\left(l^{(t)},l^{(t+1)}\right)\;,\;\textrm{otherwise}\;.
\end{array} \right.$$

%%%%%%%%%%%%%%%%%%%%%%%%%%%%%%%%%%%%%%%%%%%%%%%%%%%%%
%\begin{center}
\begin{figure}[h]
\includegraphics[scale=0.45]{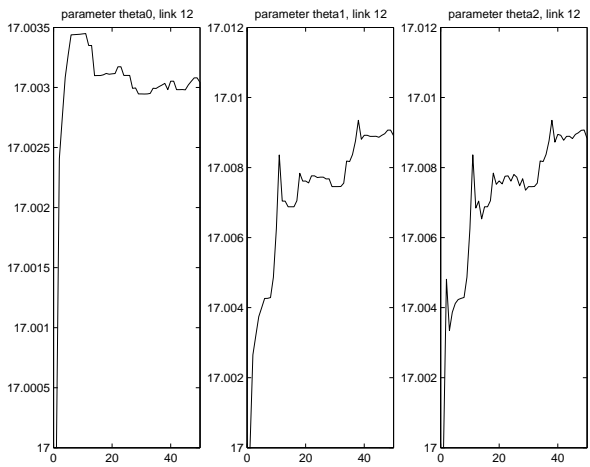}
\caption{\textbf{Parametric strategies on a MPLS-network.} We consider state spaces $\mathbf{S}^X,\;\mathbf{S}^Y,\;\mathbf{S}^Z$ of cardinality $4$. The convergence of the parameters of the strategies on the links, occures at a slower rate}
\end{figure}
%\end{center}

%\begin{center}
\begin{figure}[h]
\includegraphics[scale=0.45]{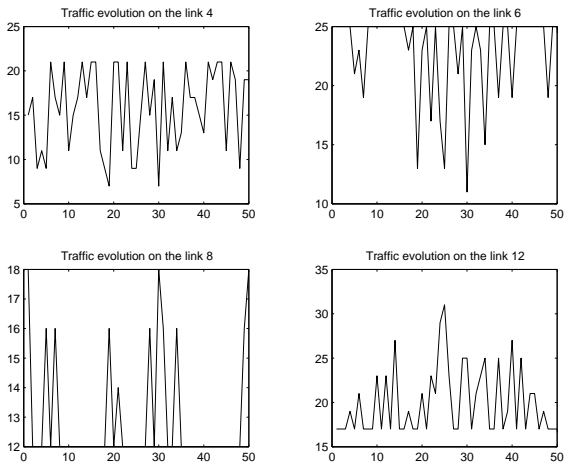}
\caption{\textbf{Parametric strategies on a MPLS-network.} Dynamic evolution of the sampled traffic trajectories on the links.}
\end{figure}
%\end{center}

%\begin{center}
\begin{figure}[h]
\includegraphics[scale=0.45]{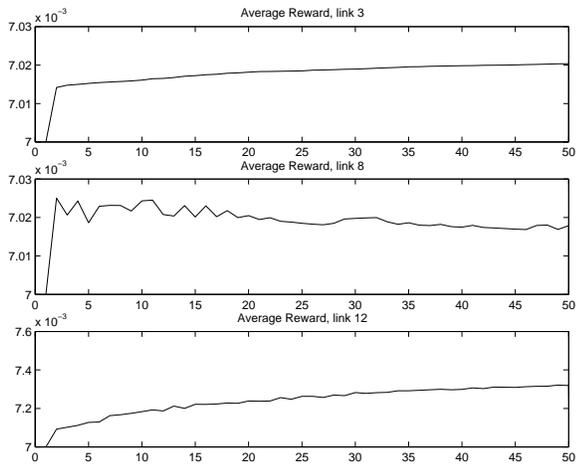}
\caption{\textbf{Parametric strategies on a MPLS-network.} Convergence of the average reward on the links.}
\end{figure}
%\end{center}

%\begin{center}
\begin{figure}[h]
\includegraphics[scale=0.45]{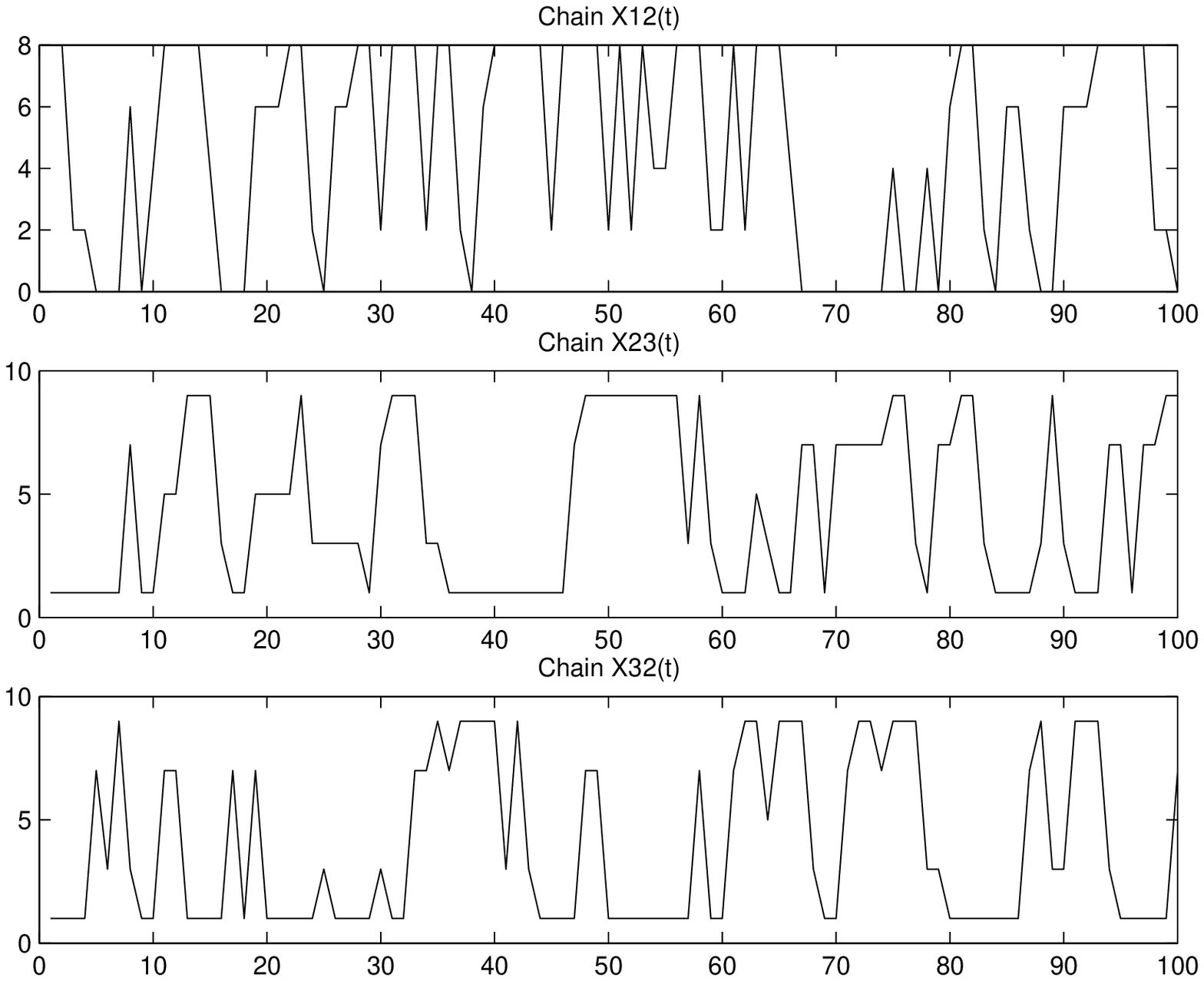}
\caption{\textbf{Parametric strategies on a MPLS-network.} Dynamic evolution of the sampled traffic trajectory on the VPN$1$.}
\end{figure}
%\end{center}
%%%%%%%%%%%%%%%%%%%%%%%%%%%%%%%%%%%%%%%%%%%%%%%%%%%%%%%%%

\bigskip

\section{Conclusion}
We have developed an original approach to tackle the problem of decision taking under uncertainty. The choice of optimizing a QoS criterion such as the delay is rather arbitrary, and could be extended to various objective functions. This article gives us rules to control optimaly a VPN so as to minimize the delay under the assumption that the traffic follows the worst possible evolution. We first determine a solution on a finite horizon $[0;T]$, using extensively Bellman's principle. But, asymptotically, we would rather apply linear programming, since under such an assumption  the strategies can be assumed stationary. The case of the management of three VPNs is also studied, via the introduction of hierarchical MDPs on the one hand, and stochastic Games on the other hand. The use of the Cross-Entropy method makes us able to forecast the trajectories of our system, provided we are given an initial state, or at least, an initial distribution on the states.\\
A curious point which could be evoked, is that the system evolved without any observation, since all the possible behaviors should be predicted and kept in memory before the system enters its initial state. In fact, the system evolution is blind, and completly disconnected of the reality. An interesting idea should be to introduce observations, so as to adapt the evolution of the system. The introduction of Partially Observed Markov Decision Processes (see $[12]$, $[13]$, $[14]$, $[15]$) might also be quite promising, but rather hard to put in application due to the large cardinality of our state spaces.\\
Indeed, the curse of dimensionality appears as soon as we have to manage a complex network of more than one VPN. The state space becomes fastly huge, and Bellman's principle gets quite difficult to put in application. Fortunately, techniques of simulation based optimization over the policy space ($[18]$), represent an alternative approach that we have tested successfully. The idea is to intoduce parametrized strategies, that depend on a set of unknown parameters. A simulation algorithm is then proposed for optimizing the average reward, and at the same time, the unknown parameters. In a practical point of view, the use of this apporach is all the more interesting, since to our knowledge, it has been tested only on few concrete case studies.   
 
\section{Bibliography}
$[1]$ DUFFIELD N.G., GOYAL P., GREENBERG A., \textit{A Flexible model for Resource Management in Virtuel Private Networks}, Proceeding of the ACM SIGCOMM Computer Communication Review, vol. $29$, pp.$95-109$, $1999$.\\
$[2]$ NALDI M., \textit{Risk Reduction in the Hose Model for VPN Design}, EURONGI Workshop on QoS and Traffic Control, Decembre $2005$.\\
$[3]$ KULKARNI V. G., \textit{Modeling, Analysis, Design, and Control of Stochastic Systems}, Springer, $1999$.\\
$[4]$ HIRIART-URRUTY J.-B., \textit{L'Optimisation}, Presses Universitaires de France, $1996$.\\ 
$[5]$ FILAR J., VRIEZE K., \textit{Competitive Markov Decision Processes}, Springer, $1996$.\\
$[6]$ ZHANG H., LIU Y., GONG W., TOWSLEY D., \textit{On the Interaction Between Overlay Routing and Traffic Engineering (MPLS), Conflicts in Routing Games}, Massachusetts University.\\
$[7]$ BEN-AMEUR W., KERIVIN H., \textit{Routing of Uncertain Traffic Demands}, Springer Science $\&$ Business Media, Optimization and Engineering, $2005$.\\
$[8]$ KORILIS Y., LAZAR A., ORDA A., \textit{Achieving Network Optima Using Stackelberg Routing Strategies}, IEEE/ACM Transactions on Networking, Vol. $5$, $1997$.\\
$[9]$ BAYNAT B., \textit{Th\'eorie des Files d'attente}, Herm\`es, $2000$.\\
$[10]$ NILIM A., EL GHAOUI L., \textit{Algorithms for Air Traffic Fow Management under Stochastic Environments}, Berkeley.\\
$[11]$ BERTSEKAS D., \textit{Dynamic Programming}, Prentice-Hall, $1987$.\\ 
$[12]$ SMALLWOOD R. D., SONDIK E., \textit{The Optimal Control of Partially Observable Markov Processes over a Finite Horizon}, Stanford Research Report, $1971$.\\
$[13]$ MONAHAN G., \textit{A Survey of POMDPs: Theory, Models, and Algorithms}, Management Science, vol. $28$,  $1982$.\\
$[14]$ LOVEJOY W., Graduate School of Business, \textit{A survey of Algorithmic Methods for POMDPs}, Annals of Operations Research $28$.\\ 
$[15]$ PINEAU J., GORDON G., THRUN S., \textit{Point-based value iteration: An anytime algorithm for POMDPs}, Carnegie Mellon University, Pittsburgh, $2005$.\\
$[16]$ DE BOER P.-T., KROESE D., MANNOR S., RUBINSTEIN R., \textit{A Tutorial on the Cross-Entropy Method}.\\
$[17]$ TOUATI C., ALTMAN E., GALTIER J., \textit{Generalized Nash Bargaining Solution for bandwidth allocation}, Elsevier, Computer Networks, Dec. $2005$.\\
$[18]$ MARBACH P., TSITSIKLIS J. N., \textit{Simulation-Based Optimization of Markov Reward Processes}, IEEE Transactions on Automatic Control, Vol. $46$, Fev. $2001$.\\
$[19]$  MARBACH P., TSITSIKLIS J. N., \textit{Approximate Gradient Methods in Policy-Space Optimization of Markov Reward Process}, Journal of Discrete Event Dynamical Systems, Vol. $13$, $2003$.\\
$[20]$ KALL P., WALLACE W. S., \textit{Stochastic Programming}, Wiley, Chichester, $1994$.\\ 

\end{document}